\documentclass{article}
\usepackage{xcolor}
\usepackage{amsmath}
\usepackage{amssymb}
\usepackage{graphicx}
\usepackage{comment}
\usepackage{setspace}
\usepackage{soul}
\usepackage{pgf}
\usepackage{tikz}
\usepackage{hyperref,mathrsfs}
\usepackage{arydshln}
\usetikzlibrary{shapes}

\usepackage{epstopdf}
\usepackage{subfigure}

\usepackage[numbers,sort&compress]{natbib}
\bibpunct[, ]{[}{]}{,}{n}{,}{,}

\newtheorem{theorem}{Theorem}[section]
\newtheorem{lemma}[theorem]{Lemma}

\newtheorem{proposition}[theorem]{Proposition}

\newtheorem{example}[theorem]{Example}
\newtheorem{remark}[theorem]{Remark}

\newenvironment{myquote}{\list{}{\leftmargin=0.3in\rightmargin=0in}\item[]}{\endlist}
\def\II{{I\!I}}
\def\lsup#1#2{\,{^{#1}}\!#2}
\def\Ic{{\overline{I}\hspace{2pt}}}
\newcommand{\ct}[1]{#1^{\mbox{\tiny\sf T}}}            
\newcommand{\cj}[1]{#1^{\mbox{\tiny\sf G}}}            
\newcommand{\N}{{\mathbb N}}                           
\newcommand{\BC}{{\mathbb C}}                          
\newcommand{\BH}{{\mathbb H}}                          
\newcommand{\BD}{{\mathbb D}}                          
\newcommand{\cO}{\mathscr{O}}                          

\def\I#1{1\hspace{-3pt} {\rm I}_{#1}}                  
\def\wek#1{\textsf{$\textbf{#1}$}}                     
\def\sp{\wek p}
\def\sq{\wek q}

\def\sd{{\wek d}}

\def\sce{\wek c}
\def\si{{\wek 1}}
\def\se{\wek e}
\def\so{{\wek 0}}
\def\sd{\wek d}

\def\su{\wek u}
\def\sv{\wek v}
\def\sg{\wek g}

\def\sr{\wek r}
\def\sR{\wek R}

\def\sN{\wek N}

\def\triu{^{\mbox{\tiny\sf U}}}
\def\bmat{\begin{pmatrix}}
\def\emat{\end{pmatrix}}

\def\ddmat#1#2#3#4{                                       
\begin{pmatrix}
#1&#2\\
#3&#4
\end{pmatrix} }
\def\dlog#1#2{\left(\genfrac{}{}{0pt}{}{#1}{#2}\right]}   
\def\enuproof#1{\begin{myquote}{\rm #1}\end{myquote}}
\def\newl{\newline\rule{0pt}{16pt}}
\def\bm#1{\mbox{\boldmath $#1$}}
\newcommand{\la}{\langle}
\newcommand{\ra}{\rangle}
\newcommand{\be}{\begin{equation}}
\newcommand{\ee}{\end{equation}}
\def\set#1{\left\{\,#1\,\right\} }                   
\newcommand{\df}{\,\stackrel{\mbox{\footnotesize\sf def}}{=}\,}
\newcommand{\imp}{\,\Rightarrow\,}
\newcommand{\tss}[1]{\textsuperscript{#1}}
\def\ca#1{{\cal #1}}
\def\Proof{\hspace{5pt}  {\bf Proof.~~}}
\def\QED{\hfill \rule{1.6mm}{1.6mm}}
\newcommand{\vspandex}{\rule[-18pt]{0pt}{36pt}}      
\newcommand{\vv}{

}                                                    

 \def\ov#1{\overline{#1}\hspace{2pt}}
\newcommand{\E}{\mbox{\sf E}\,}                           
\newcommand{\Sum}{\displaystyle\sum }
\newcommand{\Frac}{\displaystyle\frac}
\newcommand{\refrm}[1]{{\rm(\ref{#1})}}              
\newcommand{\tr}{{\rm tr}\,}                         
\begin{document}

\begin{center}
{\large\bf On signed generators of groups and algebras}
\vv

Jerzy Szulga\footnote{Department of Mathematics and Statistics, Auburn University, Alabama, USA. Email: szulgje@auburn.edu}
\end{center}
\vv\vv

\begin{quote}
{\bf Abstract}

Operators acting on the discrete random chaos yield signed multiplicative systems, extending the notion of spin matrices and quaternions. We investigate signed groups through the associated sign matrices, focusing on generators and their replacements. Of particular interest are anticommutative generators leading to a complete classification of the generated groups. The classification of finite groups of mixed commutativity is also obtained.
\vv
\vv
\vv

\end{quote}
\vv
\section{Introduction}
Studies of anticommuting matrices, with additional properties such as antisymmetry, orthogonality, roots of identity, etc., can be traced to \cite{Hur}, if not to older times, and have been attracted wide interest among both mathematicians and physicists. Even most recent works still refer to Hurwitz-Radon theorem and their consequences as well as to new approaches (cf.\ e.g., \cite{Hru}). Anticommutativity $AB=-BA$ is equivalent to the `Pythagorean formula'
\be\label{Pyth}
(A+B)^2=A^2+B^2,
\ee
and hence it suggests a specific notion of orthogonality. For example, if $A$ and $B$ are complex matrices, and one of them (say, $B$) is orthogonal, then \refrm{Pyth} implies that $A$ and $B$ are orthogonal in the classical sense ($AB=-BA\imp B^*A=-AB^*\imp \tr (AB^*)=0$). Otherwise, the  new type of orthogonality, imposed by \refrm{Pyth} (or anticommutativity),  may defy intuition (see Example \ref{2x2} below).

These ideas also emerge in the following context. The Hilbert space $L^2[0,1]$ admits a representation by Walsh series which are products of Rademacher functions and form an orthonormal basis. Analytic aspects lead to discrete Fourier analysis (cf. \cite{SWS}) while connections to quantum probability and physics entail so called `toy Fock spaces'', introduced by K.R.\ Parthasarathy \cite{Par} and further investigated by P-A.\ Meyer \cite{Mey} (for more recent developments see \cite{AttNech}). Three operators on the toy Fock space: conservation or number, creation, and annihilation, have become the  core of the theory. Yet, their commutative structure is very complicated,  although the involved calculus reveals a great deal of orderly behavior (cf.\ \cite{Mey,Par,Szu}). On the other hand, just two simple symmetries, akin to rotations, exhibit a simple commutativity and anticommutativity pattern and span all possible compositions of the quantum operators. And, yes, they were also considered by Hurwitz. This particular feature led to the notion of a signed multiplicative system \cite{Szu} which captures the simplest aspects of these and also previously considered objects and actions.   On the grounds of the group theory, the definition follows.\vv

A {\bf\em signed group}  $\ca S$, in addition to the unit 1, contains an element denoted by $-1$, commuting with all other elements; the group's members either commute or anticommute; and every member $e$ admits a {\bf signature} $e^2\in\set{-1,1}$.  Subsets of $\ca G$ or sequences of elements of $\ca G$ will be also called signed. Members $e\in\ca G$  are called {\bf\em positive} or {\bf\em negative} according to their signatures $e^2$, and the attribute extends to sets or sequences. We call a set or sequence {\bf\em pure}, if it is either positive or negative, otherwise we call it {\bf\em mixed}. \vv

Our intention is to examine such simple axiomatic system on its own, without referring to deeper and more complicated theories that may lie behind it, as far as it is possible.
Consider finitary 0-1 sequences $\wek p=(p_n)$, i.e., $p_n=0$ eventually.  A signed sequence $\wek e=(e_n)$ generates the subgroup $\ca G(\wek e)$ of powers
\[
\wek e^{\wek p}=\prod_n e_n^{p_n} \quad \mbox{(by convention, $e^0=1$ and $e^1=e$)}.
\]
A set $F\subset \ca G$ is called {\bf\em basic} and its elements are said to be {\bf\em independent} if no product of the set's elements yields $\pm 1$. We will always assume that generators $\wek e$ spanning $\ca G(\wek e)$ are basic. To reiterate: `generator' = `basic generator'. Then each member of $\ca G(\wek e)$ appears as the unique power. A generator of size $n$ entails the group of order $2^{n+1}$. Conversely, every finite signed group of order greater than 2 admits a maximal generator, and therefore the group must have order $2^{n+1}$ for some $n\ge 0$. The case $n=0$ corresponds to $\ca G_0=\set{\pm 1}$. A generator is called a {\bf\em replacement} of another generator if they generate the same group. The signatures and commutativity may vary among replacements. We identify generators that can be replaced trivially by change of signs or by permutations. For example, we identify the sequences $(\pm e_1,\pm e_2,\pm e_3)$ and the set $\set{e_1,e_2,e_3}$.  The reason is that the main use of such generators  (e.g., consisting of functions, matrices, operators, etc.) is to span the real or complex vector algebra of polynomials
\[
\set{\sum_{\wek p} a_{\wek p} \,\wek e^p: a_{\wek p}=0 \mbox{ for all but finitely many indices }}
\]
(such  vector algebra also will be called signed), where signs or permutations are irrelevant. Let us emphasize that signs should not be confused with signatures.
\vv

In Section \ref{sect:orig} we standardize the overall notation, and list necessary facts about signed groups. We indicate connections to Pauli spin matrices, quaternions, and Clifford algebras.\vv

Section \ref{sect:repl} contains  a basic dyadic linear algebra that is needed to analyze the distribution of signatures and signs reflecting commutativity and anticommutativity. The invariance of anticommutativity under replacements is examined, leading to the notion of dyadic orthogonality. The presence of at least two anticommuting elements in a signed group entails  possibly complex sign arrangements. Yet, desirably, one generator might be replaced by a `better' one that would ensure transparency and order.  \vv

In Section \ref{sect:signs} we focus on an anticommutative generator and classify groups with respect to the distribution of signatures. Up to an isomorphism and/or replacement, there are three types of odd sized generators, two types of even sized generators, and only one of an infinite (countable) generator.\vv

In the first part of Section \ref{sect:parts} we set aside the signatures in a generator as they have no effect on the commutativity matrix and find replacements and their partitions that yield orderly patterns of the generated group. We show that every finite signed group is  {\bf\em proper}, that is, it admits a generator $K\cup M$, where $K$ and $M$ commute, $K$ is anticommutative, and M is commutative. Additional patterns are examined, in particular, involving a partition of the anticommutating generator $K$ into a union of two types of doubletons, Pauli-like and quaternion-like. We classify all finite signed groups, also with respect to signatures.

\section{Notation and basic facts}\label{sect:orig}
Sequences or vectors will be denoted by a bold sanserif font. As long as the operations on elements are well defined, all operations between sequences are understood componentwise,  e.g., $\sp\sq=(p_nq_n)$. The only exception is the power $\se^{\sp}=e_1^{p_1}e_2^{p_2}\cdots$. For the sequence of powers  we will use the left superscript ${^{\wek p}\wek e}=(e_1^{p_1},e_2^{p_2},\cdots)$.  If we can add components, then we write $p=\la\sp\ra=\sum_n p_n$.\vv

 Consider the set $\BD=\set{\sp=(p_n)\in\set{0,1}^{\N}: d_n=0 \mbox{ eventually}}$ of finitary 0-1 sequences. The set is equipped with  the linear lexicographic order that coincides with the natural order of $\N$ through the binary representation $\sp\leftrightarrow n=\sum_{j=1}^\infty p_j \,2^{j-1}$, as well as with the partial componentwise order that coincides with the inclusion relation between finite subsets of $\N$, $K\leftrightarrow \sp$ with $p_n=\I K(n)$.
 \vv
  We often toggle between the multiplicative and additive notation:
\be\label{CD}
\set{-1,1}\ni c\leftrightarrow d=\frac{1-c}{2}\in\set{0,1},\quad\mbox{and also } c=(-1)^{d}.
\ee
(The function $c\to d$ may be perceived as the main branch of `the logarithm to base $-1$' on $\set{\pm 1}$.)
Additional notation will be introduced as needed, typically at the beginning of each section, and a local notation may appear right before the corresponding statement.\vv

Let us  explain how the concept emerged  (cf.\ \cite{Szu}). The real Hilbert space $\BH=L^2[0,1]$, equipped with the Walsh orthonormal basis, has been often viewed in the literature as the closure of `discrete random chaos', a linearly dense commutative algebra spanned by random signs. The random signs can be modeled by Rademacher functions that stem from the restriction of the condensed square wave, i.e.,  the periodic odd extension $r(t)$ of the indicator function $\I{[0,1)}$:
\[
r_n(t)=r(2^nt)\Big|_{[0,1]}.
\]
A finitary 0-1 sequence  $\sp$ and the Rademacher sequence $\sr=(r_n)$ yield the corresponding Walsh functions
\[
w(\sp)=\sr^{\sp}=r_1^{p_1}r_2^{p_2}\cdots\quad \longleftrightarrow \quad w_n.
\]
Bounded or even unbounded linear operators can be defined directly on the basis $\sr^{\sp}$. P-A. Meyer \cite{Mey} introduced and investigated `quantum operators', calling the Rademacher chaos a `toy Fock space'. Meyer's original notation was set-theoretic rather than algebraic.   The {\bf\em number} or {\bf\em conservation} operator $ {{N}}_j\,\sr^{\sp} =p_j \sr^{\sp}$
records the occurrences of the variable, the {\bf\em annihilation} operator  $D_j\,\sr^{\sp}     =p_j\,r_j \sr^{\sp}$ removes the variable when present in a term or the entire term lacking that variable,
the {\bf \em creation} operator ${^1\!D_j}\,\sr^{\sp} =(1-p_j)\,r_j\,\sr^{\sp}$ adds the variable when missed or puts aside the term containing the variable.

The generated algebra is isomorphic to $B(\BC^2)$, associating our operators via the usual symbols of quantum mechanics (cf. \cite{Par}):
\[
 b^- =\ddmat 0 1 0 0             \leftrightarrow D,  \quad
 b^+ = \ddmat 0 0 1 0            \leftrightarrow {^1\!D}, \quad
 n=b^o=\ddmat 0 0 0 1             \leftrightarrow {N}.
 \]
 (the subscript $j$ is fixed and then suppressed). The composition structure of the sequences $\sN,\, \wek D,\, {^1}\wek D$ is rather complex (they are just members of the standard yet sometimes unfriendly basis). In contrast, the simple symmetries,
 \[
 R_j\wek r^{\wek p} =r_j \wek r^{\wek p} \leftrightarrow \ddmat 0 1 1 0,\quad
 S_j=S_j\,\sr^{\sp}= (1-2p_j)  \,\sr^{\sp}\leftrightarrow \ddmat 1 0 0 {-1}
 \]
 together with the identity and their product  $A=RS=\ddmat 0 {-1} 1 0$  form a real anticommutative basis of $B(\BC^2)$, although with mixed signatures: $R^2=S^2=1,\,A^2=-1$.
 One can recognize the first and the third of the three Pauli spin matrices (they were already present and put to work in Hurwitz' paper \cite{Hur}).
 The products $\pm\wek R^{\wek s }\wek S^{\wek p}$ form a signed group, as defined in Introduction:
\be\label{action}
\sR^{\wek s}\wek S^{\wek p}=(-1)^{\wek s\wek p}\wek S^{\wek p}\sR^{\wek s}
\quad\imp\quad
\Big(\sR^{\wek s}\wek S^\sp\Big)\,\Big(\sR^{\wek t}\wek S^\sq\Big)=(-1)^{\wek p\wek t}\,\sR^{\wek s+\wek t}\wek S^{\sp+\sq}
\ee
In \cite{Szu}, after suppressing letters and leaving plain indices, the assignment $\dlog 0 1 = R,\, \dlog 1 0=S$ entailed a signed group of double 0-1 finitary sequences $ \dlog {\wek s} {\wek p}=\wek R^{\wek s}\wek S^{\wek p}$  that was called a `double logic'. This group is universal in the sense that it possesses a copy of every finite or countable signed system. As the original operators $R$ and $S$ can be viewed as simple $2\times 2$ matrices entailing the asymmetry $A=RS=-SR=-\dlog 1 1$, so they generate a doubling (or rather quadrupling) procedure. That is, the repetitive algorithm generates the quadruple offspring out of an element  $e$:
\[
\dlog 0 0 e=\bmat e &0\\ 0& e\emat,\quad
\dlog 0 1 e=\bmat 0 &e\\ e& 0\emat,\quad
\dlog 1 0 e=\bmat e &0\\ 0&-e\emat,\quad
\dlog 1 1 e=\bmat 0 &e\\-e&0 \emat.
\]
This yields the representation of the double sequence $\dlog {\wek s} {\wek p}$ of length $n$ as the $2^n\times 2^n$ matrix that corresponds to the iteration starting with  $e=[1]$. Although such matrices seem to be incompatible because of different dimensions for different $n$ yet the described doubling of matrices, or equivalently, concatenation of double sequences, enables one to bring matrices of different dyadic dimensions up to par. For a fixed length $n$, the simple count proves that the linearly independent matrices $\dlog {\wek s} {\wek p}$ form a basis in the vector space of matrices of size $2^n\times 2^n$.

 The matrix representation  includes generalizations of Pauli spin matrices and quaternions. For example, two components are needed to `complexify' the group. That is, focusing on sequences of length 2,
the negative element $\dlog 1 1=\dlog  {01}{01}=S_1R_1$ serves as `the imaginary number' $\imath$, i.e., the unique element besides $\pm 1=\pm\dlog {00}{00}$ that commutes with the Pauli basis
$\sigma_1=\dlog{00}{10}=R_2,\, \sigma_2=-\dlog {11}{11}=R_1S_1S_2R_2,\, \sigma_3=\dlog {10}{00}=S_2$, or the quaternion basis
${\bm i}=\imath \sigma_3=\dlog {11}{01}=S_1R_1S_2 , \,  {\bm j}=\imath \sigma_2=\dlog {10}{10}=S_2R_2, \, {\bm k}=\imath\sigma_1=\dlog {01}{11}=S_1R_1R_2$.

Note that this specific matrix representations of quaternions, one of many possible, as well as their brief `double logic' codes are the consequence of the assumed doubling algorithm. As seen above, in contrast to doubly coded implicit compositions, the explicit compositions lack transparency. Yet, the direct action on the algebraic span of the first two Rademacher functions $r_1,r_2$ (or any pair) is clear. The full matrix representations with respect to the basis $(1,r_1,r_2,r_1r_2)$ are, as expected,
\[
{\bm i}\leftrightarrow
\bmat
\cdot&1    &\cdot&\cdot\\
-1   &\cdot&\cdot&\cdot\\
\cdot&\cdot&\cdot&-1   \\
\cdot&\cdot&1    &\cdot\\
\emat,\quad
{\bm j}\leftrightarrow
\bmat
\cdot&\cdot&1    &\cdot\\
\cdot&\cdot&\cdot&1    \\
 -1  &\cdot&\cdot&\cdot\\
\cdot&-1   &\cdot&\cdot\\
\emat,\quad
{\bm k}\leftrightarrow
\bmat
\cdot&\cdot&\cdot&1    \\
\cdot&\cdot&-1   &\cdot\\
\cdot& 1   &\cdot&\cdot\\
-1   &\cdot&\cdot&\cdot\\
\emat,
\]
In addition, for doubles $D,D'\in\set{\dlog 0 0 , \dlog 0 1, \dlog 1 0, \dlog 1 1}$, we have
\[
D e \cdot D' f = (DD') ef, \quad D e\cdot Df =\dlog 0 0 ef.
\]
In particular, doubles of the same kind, or if one is $\dlog 0 0$,  preserve commutativity or anticommutativity  while the doubles of distinct kinds, different from $\dlog 0 0$, switch it. To produce more examples of anticommutating matrices we need to start somewhere.
\begin{example}\label{2x2}\rm
We will illustrate a possible advantage of using a basis alternative to the standard, Pauli, or quaternion bases.

 Let us solve the problem (cf.\ \cite{Mor}): given a $2\times 2$ complex matrix $M$ (or with entries from a commutative ring) find all anticommuting matrices $X$, i.e., $MX+XM=0$. Equivalently, the `Pythagorean formula' \refrm{Pyth} holds. To this end, let us choose the real basis $\set{I,R,S,A=RS}$, so we can easily pass to real matrices (which would be little tedious if we selected the Pauli or quaternion basis).
 We will see that the diagonal matrix $G={\rm diag}(1,1,-1)$ of signatures determines a notion of orthogonality in $\BC^3$ with the help of the bilinear form
 \[
 G(\wek v,\wek w)=\ct{\wek v}G\wek w=ax+by-cz,\quad \wek v=\ct{( a,b,c)}, \,\wek w=\ct{( x,y,z)}\in\BC^3.
 \]
 We write $M=\gamma+\wek v$, where $\wek v=aR+bS+cA$. This yields the G-conjugate or G-transpose $\cj{M}=\gamma-\wek v$. For complex coefficients we  check that
\[
(\gamma +aR+bS+cA)^2=\left(\gamma^2+a^2+b^2-c^2\right) +2\gamma (aR+bS+cA),
\]
or $(\gamma+\wek v)^2= \gamma^2+G(\wek v,\wek v)+2\gamma \wek v$. In particular, we obtain a quadratic form on $\BC^4$:
\[
M\cj M=(\gamma+\wek v)(\gamma-\wek v)=\gamma^2-G(\wek v,\wek v)=\gamma^2-a^2-b^2+c^2.
\]
Also, with the help of little algebra the `Pythagorean  formula' \refrm{Pyth} reads
\[
\gamma t+ax+by-cz=0,\quad
\gamma x=-ta,\quad \gamma y=-tb,\quad \gamma z=-tc.
\]
Thus we arrive at two cases of anticommuting matrices:
\[
\begin{array} {rll}
\mbox{(1)}&\mbox{scalar-free }\wek v \,\mbox{and}\, \wek w,&\mbox{where}\,  G(\wek v,\wek w)=0,\\
\mbox{(2)}&M=\gamma+\wek v \,\mbox{and} \, \cj M, &\mbox{where}\, M\cj M=0,\mbox{ i.e., } G(\wek v,\wek v)=a^2+b^2-c^2=\gamma^2.
\end{array}
\]
 It is easy to switch to the standard basis:
\[
\gamma+aR+bS+cA=\bmat p&q\\r&s\emat \quad\iff\quad
\begin{array}{ll}
\gamma+b=p,&\gamma-b=s,\\
     a-c=q,&a+c=r.
\end{array}
\]
Hence the restriction in Case (2) is equivalent to the singularity of the matrix, $ps=qr$. The case of commuting matrices is even simpler: a matrix $M=\gamma +\wek v$ commutes only with matrices $X=\gamma'+\alpha \wek v$.\QED
\end{example}

\vv
Both groups of matrices and of operators extend to algebras, akin to Clifford algebras,  and can be placed within the category of rigged Hilbert spaces \cite{BelTra,Szu}. In particular, any generator yields orthonormal powers and thus powers become an algebraic basis of the spanned algebra. We are not aware of simple proofs of the latter property, beyond generators of small size, without referring to the representation.
\vv
Staying merely with the definition of a signed group or its generator,  a pure anticommuting doubleton $\set{e_1, e_2}$ will be called either Pauli-like if it is positive, or quaternion-like  if it is negative. Besides this semantic link to the physical framework,  from now on we will not use any matrix or operator representation of signed groups, only the basic definition.
\section{Replacements}\label{sect:repl}

\subsection{The AC (anti/commutativity) matrix}
The signatures of independent  elements  define the sequence $\bm\sigma=\sigma(\se)$  and the symmetric commutativity function on $\ca S\times \ca S$ defines the symmetric matrix $C=[c(e,e')]$, often called the {\bf AC-matrix}:
 \[
c(e,e')=e\,\circ\, e'=\left\{
\begin{array}{rl}
1,&\mbox{if~ $e\,e'=e'\,e$},\\
-1,&\mbox{if~ $e\,e'=-e'\,e$}.\\
\end{array}
\right.
\]

The sequence $\bm\sigma$ has no effect  on the AC-matrix but the former depends on the latter.
 For a square matrix $C=[c_{jk}]$, denote by $C\triu=[c_{jk}\,\I{\set{k>j}}]$ the upper triangular cut of $C$ (Matlab's `{\tt triu(C,1)}').
\begin{lemma}\label{swaps}
Let  $\wek e$ be a basic sequence in $\ca S$ with an AC-matrix $C=[c_{jk}]=[c(e_j,e_k)]$ or its 0-1 equivalent $D$ (cf.\ \refrm{CD}).
Then
\[
\wek e^{\wek p}\,\wek e^{\wek q}=(-1)^{\langle {D\triu\sp}, {\sq}\rangle}\,\wek e^{\wek p+\wek q} , \quad \sp,\,\sq\in\BD,
\]
that is,
\[
{\langle {D\triu\sp}, {\sq}\rangle}=\sum_{j\ge 1}\,\sum_{k>j}d_{jk}{p}_k\,{q}_j=\frac{1}{2}\sum_{j\ge 1}\,\sum_{k>j}\left(1-{c}_{jk}\right){p}_k\,{q}_j.
\]
In particular, we observe the following properties:
\begin{enumerate}
\item $
\se^\sp\,\circ\,\se^\sq=(-1)^{(\la C\sp,\sq\ra-\sp\sq)/2};
$
\item for an anticommutative sequence,
\be\label{signpq}
\se^\sp\,\circ\,\se^\sq=(-1)^{\la\sp\ra\la\sq\ra-\sp\sq};
\ee
\item
The only group $\ca G\subset \ca S$ with entirely anticommutative  subset $\ca G\setminus\ca G_0$ must be of order 8, i.e., must be generated by just two elements (in other words,  it must be isomorphic either to the Pauli group or to the quaternion group);
\item\label{diag4}
Denote $p=\la \sp\ra$. The signature $\sigma(e^{\sp})=\se^{2\sp} (-1)^{p(p-1)/2}$. That is, for pure sequences, \begin{enumerate}
\item
If $\se$ is positive then $\sigma(\se^\sp)=-1$ iff $p=2 \,(\mbox{\rm mod }4)$ or $p=3\, (\mbox{\rm mod } 4)$;
\item
If $\se$ is negative then $\sigma(\se^\sp)=-1$ iff $p=2\,(\mbox{\rm mod } 4)$ or $p=1 \,(\mbox{\rm mod } 4)$.
\end{enumerate}
\item
Let $s_+$ and $s_-$ denote the counts of negative elements in  $\set{\wek e^{\wek p}:\sp\in \BD_n}$ when $\se$ is positive or negative, respectively. Then
\be\label{ss}
\mbox{$s_+=b_2+b_3$ and $s_-=b_1+b_2$}, \, \mbox{where}\,  b_q=\Sum_{j} {n\choose q+4j}\,\mbox{ (i.e., $q+4j\le n$)}.
\ee
\end{enumerate}
\end{lemma}
\Proof
First, we note the tautology for $i,j\in \ca S$ and $p,q\in\set {0,1}$:
\[
 i^{p}\,j^{q}=
 \left\{
 \begin{array}{ll}
i^{p}\,i^{q}=i^{{p}+{q}},&\mbox{if $ \,i=j$};\\
 c(i,j)\,j^{q}
\, i^{p},&\mbox{if $ \,i\neq j$}.\\
\end{array}\right.
\]
Consider
\[
e_1^{{p}_1}\cdots e_n^{{p}_n}\,\,\,e_1^{{q}_1}\cdots e_n^{{q}_n}
\]
While moving the factor $e_1^{{q}_1}$ to the left until it meets $e_1^{p_1}$, the tautology records consecutive swaps and produces the exponent $(d_{2,1}\,p_2+\cdots+d_{n,1}\,p_n)q_1$. The repetition of the process for each next factor entails the quadratic form in the exponent of $(-1)$. The remaining formulas follow immediately. Particular forms of the exponents are due to the identities $(-1)^a=(-1)^{-a}$ and $(-1)^\sp=(-1)^{\la\sp\ra}$.
\vv
(3) Three  anticommuting independent elements $e_1, e_2, e_3$ would  yield $e_1(e_2e_3)=(e_2e_3)e_1$.\vv

(4) The classification simply relies on the parity of the quadratic functions $p(p-1)/2$ and $p(p+1)/2$.\vv

(5) The identities follow directly from statement (4).\QED
\vv
\begin{remark}\rm{~}

\begin{enumerate}
\item The exponents $\sp$ and $\sq$ correspond uniquely to the finite subsets of $\N$, and hence for an anticommutative sequence the sign in \refrm{signpq} is determined by the parity of the off-diagonal elements in the Cartesian product:
\[
\I P=\sum_j p_j\I{\{j\}},\quad \I Q=\sum_j q_j\I{\{j\}}\quad\imp\quad \la\sp\ra\la\sq\ra-\sp\sq=|P\times Q|-|P\cap Q|.
\]
\item The order restriction, addressed in Statement 3, is  a consequence of the imposed condition: anticommutitivity of all basic elements of the induced algebra. In particular, there appears a quick proof of an analogous statement \cite[Prop.\ 4.5]{JadSzu}. However, the formal simplicity, bordering on triviality, comes at a cost: the argument for Statement 3  is devoid of a deeper physical meaning. If the requirement of total anticommutativity of the induced products is relaxed then the dimension size (but not its numerical structure, though) is no longer limited.\QED
    \end{enumerate}
\end{remark}

We will take a closer look at the replacements of generators of fixed length, conforming to notational conventions of linear algebra. In the literature devoted to algebra such topics rarely appear on their own, being rather a margin of algebra of finite fields. Properties in this particular case of fields of characteristic 2 are usually derived from the general theory as much as it is possible. Often the beautiful general theory fails to deliver because its tools, e.g., quadratic forms and polynomials, break at the boundary of characteristic 2. Let us illustrate the situation by quoting Robert Wilson \cite[Chap. 3.8]{wil}:
\begin{quote}
{\em In characteristic 2 everything is different. The quadratic form has a different definition, the canonical forms are different, there are no reflections, the determinant tells us nothing, and there is no spinor norm.}
\end{quote}
Another area that deals specifically with dyadic matrices is the coding theory, equipped with its own specialized terminology  (e.g., `code' is the synonym of a dyadic vector space) and its own topics of interest. We would rather refer to the common linear algebra and basic arithmetics, keeping the simple - simple. The practice of dyadic vector spaces is like the classical practice, but often not quite the same.
\vv
Thus, first we refresh and adjust our old notation.
Our default form for a vector $\sv$ of size $n$ is a vertical $n\times 1$ matrix, and $v=\la\sv\ra=\sum_i v_i$. If $\wek v$ has 0-1 components then $\la\sv\ra$ simply counts ones, so it may be called the {\bf \em mass} of $\sv$. For a matrix $V=[v_j^i]$, $\sv_j$ denote its columns while $\sv^i$ denote its rows (horizontal vectors), entailing $\sv=\ct{[v_1,\dots, v_n]}$ and $\sv^{\boldsymbol{\cdot}}= \ct{[v^1,\dots, v^n]}$, and also
$v=\la \sv\ra=\la\sv^{\boldsymbol{\cdot}}\ra=\sum_{i,j}v_j^i$. The zero vector is denoted by $\so$ and $\si$ is the vector of ones, while $\si_j$ has all zeros except for the solitary 1 at the $j$\tss{th} position.  The zero matrix is denoted by $O$, and we denote a block of ones by  $\II\df \si\ct{\si}$ (or `{\tt ONES(n,k)}' in the Matlab's notation).  \vv

Two integers $m,\,n$ are called {\bf\em congruent},  which is written as $m\cong n$, when $m=n$ (mod 2). For $c\in\set{0,1}$, we denote its complement by $\ov c=1-c$.  The relations are extended to integer matrices.  That is, $P\cong Q$ if $P=Q+2R$ for some integer matrix $R$ and $\ov P=\II-P$. Even though we focus on dyadic matrices yet occasionally we still wish to keep their `memory' of being integer matrices.\vv

 A binary or integer matrix $P$ with columns $\sp_1,\sp_2,\dots$ transforms a  signed anticommutative basic sequence $\se$ to
the sequence $^P\se\df (\se^{\sp_1},\se^{\sp_2},\dots)$. We use the left superscript because the right superscript would yield a power (cf.\ the notation section). Note that
$^Q(^P\se) =(\pm \se^{\sq^1P},\pm\se^{\sq^2P},\dots)$, associated with but not exactly equal to $^{QP}\se$. By convention, we identify both sequences.\vv

 Formula \refrm{signpq} provides the characterizing condition for  anticommutativity of the transformed sequence. In the language of matrices the condition reads
 \be\label{V}
 \ct P\, \ov  I\, P \cong \ov  I.
 \ee
We simply call such $P$ anticommutative.  The set of anticommutative  integer matrices $P$ is closed under the product, i.e., it is a semigroup. We will examine anticommutative $P$ such that $^P\se$ replaces the original generator, preserving the independence.\vv

Clearly, anticommutativity is preserved by permutations of columns or rows of $P$, and the resulting matrix $Q$ will be deemed {\bf\em indistinguishable} from $P$, $Q\leftrightarrow P$ in short.\vv
\subsection{Dyadic inverse}
Let us note a crude multiplication table for rectangular matrices, where $n$ is the inner dimension:
\be\label{multt}
\mbox{if $n$ is even:
\,\,
\begin{tabular}{r|cc}
$\cong$ & $\Ic$ & $\II$ \\ \hline
$\Ic $  & $I$   & $\II$ \\
$\II$   & $\II$ & $O$
\end{tabular},}
\qquad
\mbox{if $n$ is odd:
\,\,
\begin{tabular}{r|cc}
$\cong$ & $\Ic$ & $\II$ \\ \hline
$\Ic $  & $\Ic$ & $O$   \\
$\II$   & $O$   & $\II$
\end{tabular}}.
\ee
Hence $\ov  I\ov  I\ov  I\cong\ov  I$ regardless of parity. Elementary computations stand behind these congruence formulas, e.g., for square $n\times n$ matrices,
\[
\II^2=n\II,\quad \ov  I^2=I+(n-2)\II,\quad \ov  I\II=\II\ov  I=(n-1)\II,\quad\mbox{etc.}
\]
Notice that $\ov  I$ is nonsingular but $\Ic^{-1}=(n-1)^{-1}\II-I$ is not an integer matrix. In view of this observation, tracing the idea of replacement, we say that $P$ is {\bf\em dyadically invertible}, DI in short, if there is an integer matrix $A$ such that $AP\cong I$.

\begin{proposition}\label{dyam}
 Let $P$ be a square integer matrix.
\begin{enumerate}
\item If $P\cong I$  then $\det P\neq 0$.
\enuproof
{\rm Indeed, using the Levi-Civita symbols $\epsilon_{\wek j}$,
 $
n!\,\det P= \Sum_{\wek j}\epsilon_{\wek j}\, p_{1j_1}\cdots p_{nj_n},
 $
 where the sum runs over $n!$ permutations of $(1,\dots, n)$. If $P\cong I$ then with the solitary exception of an odd integer $p_{11}\cdots p_{nn}$ all products are even, so the sum is odd. }
  \item If $Q$ is DI and $P\cong Q$, then $P$ is also DI.
\enuproof
    {\rm Indeed, $AQ\cong I$ and $P=Q+2R$ entails $AP=AQ +2AR\cong I$.}
 \item
The following conditions are equivalent:
 \begin{enumerate}
 \item $P$ is DI;
 \item $P\sce\cong 0\imp \sce\cong 0$;
 \item $\det P\cong 1$  and $(\det P)\cdot  P^{-1}$ is an integer matrix;
 \item There is an integer matrix $B$ such that $PB\cong I$.
 \end{enumerate}
\enuproof
 {\rm (a) $\imp$ (b): Let $AP=I+2R$. If $P\sce=2\sr$ then $\sce+2R\sce =AP\sce=2A\sr$, i.e., $\sce\cong 0$.
\newl
 (b) $\imp$ (c): (b) ensures no even column. In fact, by (2) we may consider just a binary 0-1 matrix with no zero column. Therefore a Gauss row reduction based only on permutations and subtraction of rows yields a triangular matrix with ones on the diagonal. The other formula just uses the adjugate matrix.
\newl
(c) $\imp$ (d): choose $B={\rm adj}(P)=(\det P)\,P^{-1}$ in lieu of the actual inverse.
\newl
(d) $\imp$ (a). Assume (d), i.e. $PB=I+2R$ for an integer matrix $R$. Hence $b=\det B\cong 1$ and $b\, B^{-1}$ is integer by (a) $\imp$ (c). Hence $b\,P=b\,B^{-1}+2b\,RB^{-1}$. Put $A=B$, so $AP\cong b\,BP=b\,I+2b\,BRB^{-1}\cong I+2R'$. \QED}
\end{enumerate}
\end{proposition}

\begin{remark}\label{remDI}\rm
Let $P$ be a square integer matrix.
\begin{enumerate}
\item
As integer matrices the left dyadic inverse $A$ and the right dyadic inverse $B$ may differ but they are congruent. In particular,  in defining $P$ as {\bf\em dyadically orthogonal} (or `D-orthogonal', in short) by the relation $\ct PP\cong I$,  the condition is equivalent to $P\ct P\cong I$.
\item
$P^{-1}$ is integer iff $|\det P|=1$.
\item
Let $P\cong \Ic$. If $n$ is even, then $P$ is D-orthogonal, so $P$ is DI and thus nonsingular. If $n$ is odd then $\det P\cong 0$, so $P$ is not DI and may be singular. However, the actual $\Ic$ is nonsingular.

\enuproof
{The first statement is the particular case of Remark 1 above. Consider an odd $n$. Let us mark the dimension, $P=P_n$. We do not know a quick argument, so a quite tedious row reduction followed by equally tedious column reduction will reduce $\det P_n$ to $\det P_{n-2}$, where $P_{n-2}\cong \Ic_{n-2}$. Repeating, we reduce the case  to dimension $n=3$, and then to dimension 1 with $P_1=[p]$ for some even number $p$, possibly 0. E.g., $P_3=[0,1,1;-1,0,1;-1,-1,0]\mapsto P_1=[0]$.
\newl
The inverse of the actual $\Ic$ was shown right above the proposition.}
 \item
 While every DI matrix is nonsingular, the inverse implication fails.

\enuproof
{\rm Indeed, the inverse is just a rational matrix, not necessarily an integer matrix. A counterexample may involve just even rows, for an odd $n$ we take $\Ic$ while for an even $n$, say, $n=4$:
 \[
 P=\bmat 1 & 1 & 1 & 1 \\ 1 & 1 & 0 & 0\\ 1 & 0 & 1 & 0\\ 1 & 0 & 0 & 1 \emat
\quad\imp\quad\mbox{$P\wek 1\cong \wek 0,\quad \det P\neq 0$}.
 \]}
 \item
 If $Q$ is another square integer matrix and $PQ$ is DI, then both $P$ and $Q$ are DI. In particular, if an anticommutative $P$ is DI, then its dyadic inverse is also anticommutative.

 \enuproof{Indeed, considering integer determinants: the odd product requires odd factors.}
    \item
The set $\ca D$ of DI matrices forms a transpose-invariant group.
\QED
\end{enumerate}
\end{remark}
The columns of a fixed DI matrix form a basis in the vector space of all DI matrices. We need to decide whether or not we identify DI matrices modified by permutations $P\leftrightarrow \ct \pi P\pi$, i.e., whether we perceive an ordered or undordered basis. The number of ordered bases (cf. \cite[Lemma 1.17]{car} or \cite[(3.1)]{wil}) can be found exactly as well as its asymptotics:
\[
\prod_{i=1}^{n-1}\left(2^n-2^i\right)=d_n\,2^{n^2},\quad\mbox{where }d_n=\prod_{i=1}^{n-1}\left(1-2^{-i}\right)\to \phi(1/2)\approx 0.288788095...,
\]
and $\phi(z)$ is the Euler function. That is, almost 30\% of large dyadic matrices are DI. The division by $n!$ yields the number of unordered bases or permutation equivalent DI matrices.

\vv
An anticommutative matrix admits yet another characterization:
\be\label{PTP}
C=\ct P P\cong\bmat I &\II\\ \II & \Ic \emat.
\ee
Factually,  the matrix on the right is $\ct \pi C\pi$ but we use the inverse permutations for the more transparent display, which will be called {\bf \em canonical} from now on.

\begin{example}\label{addcolrow} {~}\rm

1. Consider the matrix of the generator's transformation `{\em multiply by an element $e_i$}', say, by the 1\tss{st} element $e_1$. That is, the first column is added to all other columns. The corresponding matrix $C_1=\bmat 1 & \ct{\wek 1}\\ \wek 0 & I \emat$ is anticommutative:
\[
PC_1=[\sp_1|\sp_2+\sp_1|\dots|\sp_n+\sp_1],\quad
\ct{C_1}C_1\cong \bmat 1 & \ct{\wek 1}\\\wek 1 & \Ic \emat=\bmat 1 &\lsup 1\II_{n-1}\\ \lsup {n-1} \II_1 &\Ic\emat.
\]
However, $\ct{C_1}$ (deemed as the left action $P\mapsto C_1^T P$, adding the first row to all other rows) is anticommutative iff $n$ is even, because $C_1\ct{C_1}\cong \bmat n & \ct{\wek 1}\\\wek 1 & \Ic \emat$.\vv
Factually, this example involves the replacement
\be\label{Ff}
F(f)=\set{f}\cup f(F\setminus\set{f})=\set{f}\cup\set{ff':f'\neq f},\quad f\in F.
\ee
We observe that the transformation preserves anticommutativity.\vv

2. Consider the `{\em add one row to another}' matrix $R$ ($P\mapsto RP$), say, the $1$\tss{st} row to the $2$\tss{nd} row:
\[
R=\lsup 1 R_2= \bmat 1 & 0 &\ct{\so} \\ 1 & 1 &\ct{\so} \\ \so & \so & I\emat\quad\imp\quad
\ct R R\cong \bmat 0 & 1 &\ct{\so} \\ 1 & 1 &\ct{\so} \\ \so & \so & I\emat
\]
Since such matrices do not anticommute, the Gauss row reduction does not preserve anticommutativity.\QED
\end{example}

The identity $I$ occupies the entire matrix \refrm{PTP} iff $P$ is D-orthogonal. The other extreme occurs when $C=\ct PP=\Ic$, so we are tempted to call such $P$ `{\bf\em antiorthogonal}'.

\begin{proposition} Let $P$ be anticommutative  of dimension $n$.

Consider the matrix $C$ from \refrm{PTP} and let $n=m+k$, where $m$ is the dimension of the NW (North-West) block and $k$ is the dimension of the SE block. Then, with the only exception of even $m$ and odd $k$ (making $n$ odd), $P$ is DI.
\end{proposition}
\Proof
 We will examine $C^2$, going to $C^3$ if necessary. It will be convenient to  mark the connecting dimensions of essence by `0' (even) or `1' (odd) in the multiplication table \refrm{multt}:
\be\label{multt01}
\mbox{\begin{tabular}{r|cc}
          & ${^0}{\Ic}$ & $\lsup 0{\II}$ \\ \hline
$\Ic_0 $  & $I$            & $\II$ \\
$\II_0$   & $\II$          & $O$
\end{tabular}}
\qquad \mbox{and}\qquad
\mbox{\begin{tabular}{r|cc}
          & ${^1}{\Ic}$ & $\lsup 1 {\II}$ \\ \hline
$\Ic_1 $  & $\Ic$          & $O$   \\
$\II_1$   & $O$            & $\II$
\end{tabular}}.
\ee
Depending on the parity of $m$ and $k$, in $C^2$ we encounter the following four patterns:

\def\cases{
\begin{tabular}{r|cc}
  & 0     &   1     \\ \hline
0 & $I_0$ & $\Ic_1$ \\
1 & $I_1$ & $D$
\end{tabular}}
\def\mdown
{\begin{tabular}{r}
\rule{0pt}{10pt}\\
m:
\end{tabular}}
\def\Dright
{\begin{tabular}{l}
\rule{0pt}{10pt}\\
  where $D= \bmat \Ic_1 &\lsup 1{\II}_1 \\ \lsup 1{\II}_1 & I_1\emat$.
\end{tabular}}
\begin{center}
\!\!\!\begin{tabular}{rcl}
                 &\rule{12pt}{0pt} k: &       \\
\mdown           & \cases             &\Dright\\
\end{tabular}
\end{center}
We verify that $C^3\cong CD\cong I$. Then, the statement follows from Remark \ref{remDI}.5.\QED

 \begin{remark} Let $P$ be anticommutative with $n=m+k$ defined above.
 \begin{enumerate}
 \item
Necessarily, $m\ge 1$, i.e., there is no antiorthogonal $P$ of even dimension. Further, an antiorthogonal matrix is never DI.

\enuproof
{ The analysis of parity of $m$ admits $m=0$, which does not exclude DI when $k$ is even. However, an antiorthogonal matrix $P$ cannot have a dyadic inverse. Indeed,  suppose by contrary that $PQ\cong I$ for some integer matrix $Q$. Then  $I\cong \ct{(PQ)} (PQ)= \ct Q(\ct P P)Q\cong  \ct Q\Ic Q =\ct Q\II Q-\ct Q Q= \sq\ct{\sq} -\ct Q Q.$ This would yield the congruence of diagonals, $\si\cong \so$, a contradiction.}
 \item The complement $P\leftrightarrow \ov P$ preserves both extremes of dyadic-orthogonality when $n$ is even, and switches them if $n$ is odd.
\enuproof
 {\rm Indeed, $\ct{\ov P}\ov P\cong n\II+\ct P P$. }\QED
\end{enumerate}
 \end{remark}

 \subsection{Orthogonal replacements}
 We will examine some properties of D-orthogonal matrices, $\ct P P\cong I$, including generating algorithms. We begin by counting.\vv

Denote the group of D-orthogonal matrices by $\cO$.  Writing $P=I+A$, the D-orthogonality means exactly that
 \[
 A+\ct A +\ct A A\cong O\qquad\mbox{(for symmetric matrices: $\ct AA\cong 0$)}.
 \]
  Row and column permutations $P\leftrightarrow \pi P\sigma$ scatter the elements, preserving the orthogonality but changing the appearance. Although the result may be deemed indistinguishable from the original but the mapping is not a homomorphism in $\cO$, $P\leftrightarrow P',\,Q\leftrightarrow Q'\not\Rightarrow PQ \leftrightarrow P'Q'$, with the exception of a  fixed $\pi$ and $\sigma=\ct\pi$. \vv

  The order of $\cO$ is well known and elementary calculations can be found, e.g., in \cite{MacW}.  Two DI matrices $M$ and $N$ are called {\bf\em orthogonally equivalent}, $M\sim N$ if $NM^{-1}\in \cO$, i.e., $N\in \cO M$. There is a one-to-one correspondence between the coset $\cO M$ and the set of DI symmetric matrices $\ct MM$, so this set and $\cO$ are equal in size. Then the author's \cite{MacW} count of symmetric matrices relies on a classical 1938  result \cite{Alb} stating that the representation $\ct MM$ of a symmetric matrix is possible if and only if the diagonal of this symmetric matrix has nonzero mass. The count, including also arbitrary ranks, is performed in general context of a finite field $GF(q)$, $q=2^p$, $p$ a prime. In our case $q=2$ the number of DI symmetric matrices of dimension $n=2m$ or $n=2m+1$ equals
  \[
  N_m=\prod_{i=1}^m 2^{2i} \left(2^{2i-1}-1\right) =c_m\,2^{2m^2+m},\quad \mbox{where }c_m\to \frac{\phi(1/2)}{\phi(1/4)}\approx 0.41942244...,
  \]
  using the Euler function $\phi$ again. At the same time the number of DI symmetric matrices with 0 diagonal is $2^{-2m} N_m$ when $n$ is even and there is none when $n$ is odd. Hence, the number of D-orthogonal matrices is still in the range of $2^{2m^2}$.  Therefore, even when pooling together  permutation modified matrices (involving the division by a number not exceeding $(n!)^2$, the order does not change significantly (which is subject to discussion about the meaning of `significance', involving the degree of smallness among big numbers). Thus, we obtain a crude estimate of the order $\approx 2^{n^2/2}$.

 \subsubsection{A modified Gram-Schmidt algorithm}
 First, we will discuss the analog of the Gram-Schmidt orthogonalization which turns out to  be almost a copy of the classical one but with a twist. It might be tempting to dub two dyadic vectors $\sp$ $\sq$ D-orthogonal if $\la \sp\sq\ra=0$, or even $C$-orthogonal  with the help of the bilinear form $\ct\sp C\sq=0$, but soon that would lead to a sort of havoc, or a distraction, or an ambiguity at least. Such situations are discussed in detail in \cite[3.4-3.8]{wil}. \vv

However, here D-orthogonal matrices popped up naturally prior to defining D-orthogonal vectors.  For any dyadic matrix $C$, $\ct CC$ contains the masses of columns of $C$ on the diagonal and the masses of their intersections (or of entry-wise products) off the diagonal. So, a D-orthogonal matrix by definition has odd columns (i.e., of odd masses) with even mutual intersections. So, it seems natural to say that vectors are D-orthogonal if they are columns of a D-orthogonal matrix, i.e. $\sp {\perp} \sq$ if $\la\sp\sq\ra\cong 0$ \underline{and} $\la\sp\ra\cong \la\sq\ra\cong 1$, which entails the notion of the orthogonal complement of an odd dyadic vector. At the same time, an even vector would have no D-orthogonal complement by default. \vv

The vector $\si$, a `flatline',  whether even or odd, has no orthogonal complement, nor does any vector subspace  containing it. The classical Gram-Schmidt process seems to be useless. Even low dimension examples show that a subspace may have not even a pair of orthogonal vectors, or just a few. It may contain the `flatline' or not at all. Yet, the Gram-Schmidt process can be used to construct orthogonal matrices, i.e., orthogonal bases, starting with a single vector, which we may also call `D-orthogonal' although it is solitary.
\vv

Given D-orthogonal vectors $\sp_1,\dots,\sp_k$ of length $n$ such that
\be\label{d-span}
\si\notin {\rm span} \set{\sp_1,\dots,\sp_k}
\ee
we put $\sp=\sp_1+\cdots+\sp_k$. Necessarily, $k<n$. Next,  consider the projection of a vector $\su$:
\be\label{proju}
\E[\su|\sp_1,\dots,\sp_k]=\sum_{j=1}^k \la \su\sp_j\ra\,\sp_j \quad\mbox{and}\quad \su'=\su+\E[\su|\sp_1,\dots,\sp_k].
\ee
We check that
\[
\la\su' \sp_j\ra\cong 0\quad\mbox{and}\quad \la\su'\ra=\la\su\ra+\la\su'\sp  \ra.
\]
The vector $\su'$ is odd when exactly one of the summands is odd. Therefore and since  $\sp\neq\si$ we have at least two choices  of such $\su$ (say, odd or even).  Because of two choices we can augment our orthogonal set by $\sp_{k+1}=\su'$ to avoid the `flatline', i.e., $\sp_1+\cdots +\sp_{k+1}\neq \si$ while $k+1<n$. Then the algorithm continues until it is complete when $k+1=n$.\QED

\subsubsection{Quick generators}
To obtain a quick example of a D-orthogonal matrix consider
\be\label{simP}
P\cong \bmat I & O\\O & \Ic\emat,
\ee
 where the lower block has even dimension. Given two 0-1 sequences $\su,\sv $,  the scattered block of ones can be written with the help of the classical tensor notation $\lsup{\su}\II_{\sv}=\su\ct{\sv}=\su\otimes \sv$.
In particular, $\lsup {\sv}\Ic_{\sv}=\ov{\sv\otimes \sv}$, hence the old multiplication table \refrm{multt} now takes the following form:
\be\label{IIscatter}
\begin{array}{rcl}
\su\otimes\sv\cdot\sr\otimes \wek s         &=& \la \sv\sr\ra\, \su\otimes \wek s\\
\su\otimes\sv \cdot\ov{\sv\otimes \sv}      &=&  v\,\su\otimes \sv\\
\ov{\sv\otimes \sv}\cdot\ov{\sv\otimes \sv} &= & \left(\ov{\sv\otimes \sv}\right)^v.
\end{array}
\ee
(recall that $a^0=1$ and $a^1=a$). More general formulas are available but we will not need them here. So, $I+\su\otimes \su$, where $\su$ is even, is permutation equivalent to \refrm{simP}. Let us examine a slight extension  $P= I+\sum_{i} \su_i\otimes\sv_i$,
 where  $\la\su_i\su_j\ra\cong 0$ (even intersections, including self). First, we check when $P$ is dyadically orthogonal, $I\cong \ct P P=$
 \[
I+\sum_{i}\left(\sv_i\otimes\su_i +\su_i\otimes\sv_i\right)+\sum_i\sum_j  \la \su_i\su_j\ra\, \sv_i\otimes\sv_j\cong I+\sum_{i}\left(\sv_i\otimes\su_i +\su_i\otimes\sv_i\right).
 \]
We see  that  necessarily $\su_i\cong \sv_i$, i.e.,
 \be\label{dorts}
P\cong I+A=I+\sum_{i} \su_i\otimes\su_i.
\ee
When the vectors $\su_i$ are factually disjoint, a permutation equivalent form of $A$ will have quadratic blocks of ones along the diagonal. However, in general it is harder to visualize its form even with the help of permutations.\vv

The family of D-orthogonal matrices \refrm{dorts}, denoted by $\ca K$,
is not closed under the  product $PQ$ (first of all, the symmetry is not product invariant):
\[
\begin{array}{c}
\vspandex
\left(I+\Sum_{i} \su_i\otimes\su_i\right)\left( I+\Sum_{i} \sv_i\otimes\sv_i\right)\\
=I+ \Sum_{i} \left(\su_i\otimes\su_i +\sv_i\otimes\sv_i\right)+ \Sum_{i} \Sum_j \la \su_i\sv_j\ra\, \su_i\otimes\sv_j\\
=I+ \Sum_{i} \left(\su_i\otimes(\su_i+\wek t_i) +\sv_i\otimes\sv_i\right),\quad
\mbox{where}\quad
\wek t_i= \Sum_j \la \su_i\sv_j\ra\, \sv_j.
\end{array}
\]
For example, $P=I+A,\,Q=I+B$, where
\[
A=\bmat 1&1&0\\ 1&1&0\\0&0&0\\\emat,\quad  B=\bmat 0&0&0\\ 0&1&1\\0&1&1\\\emat\quad
\imp\quad A+B+AB=\bmat 1&0&1\\ 1&1&0\\0&1&1\\\emat
\]
However, here $P, Q$ here are permutations, and their product is a permutation, all indistinguishable from $I$.
A simplifying assumption, e.g., the even intersections $\la \su_i\sv_j\ra \cong 0$, yields the same type:
\be\label{dorts2}
PQ\cong I+ \Sum_{i}\su_i\otimes\su_i +\Sum_i \sv_i\otimes\sv_i=I+\Sum_i \su'_i\otimes\su'_i.
\ee
So,  $\ca K^2$ is distinguishable from $\ca K$, e.g., for  odd $n$ and  odd blocks, we may obtain an asymmetric matrix:
\[
\bmat \Ic & \II & O\\ \II & \Ic & O \\ O & O & I\emat\cdot
\bmat I & O & O\\ O  & \Ic & \II \\ O & \II & \Ic\emat=
\bmat \Ic & O & \II\\ \II  & \Ic & O \\ O & \II & \Ic\emat\cong
I+\bmat \II & O & \II\\ \II  & \II & O \\ O & \II & \II\emat\notin \ca K
\]
When all blocks are even, as expected in \refrm{dorts2}, the product equals
\[
\bmat \Ic & \II & O\\ \II  & I & \II \\ O & \II & \Ic\emat\cong
I+\bmat \II & \II & O\\ \II  & O & \II \\ O & \II & \II\emat\in \ca K.
\]
Although $\cO$ is large, it is finite, so the sequence $\ca K^k$ must eventually become constant. Two questions emerge, the first one is whether the entire $\cO$ can be reached by this procedure. If the answer is affirmative, the second question is how fast this task can be achieved. However, we observe that the number of simple matrices \refrm{dorts} is just a humble big number, compared to the total count.
 \begin{proposition}
 Denote by $p_0(n)$ the number of distinguishable $P=I+A$ of type \refrm{dorts} with disjoint vectors. Then the range of $p_0(n)$ is approximately equal to $c\, e^{d\sqrt{n}}$ as $n\to \infty$, for suitable constants $c,d>0$.
 \end{proposition}
 \Proof
 The number of canonical $P$s is related to the partition function $p(m)$, the number of ways $m$ can be written as a sum of positive integers with order irrelevant (cf. \cite[(23)]{wolfram} for the definition and properties).  When $n$ is odd, $n=2m+1$, we must select an odd number $2j+1$, $0\le j\le m$, to make the dimension of $I$. When $n$ is even, $n=2m$, we may either  select $I$ of dimension $2j$, $j\le 2m$, or skip it. Therefore, the number $p_0(n)$ equals
 \[
 \begin{array}{rl}
 \mbox{if $n=2m+1$:} & p_0 (n)=\Sum_{j=0}^m j\,p(m-j),\\
 \mbox{if $n=2m$}:   & p_0 (n)=\Sum_{j=0}^m j\,p(m-j)+p(m).\\
 \end{array}
 \]
 The Hardy and Ramanujan (1918) approximation \cite[(23)]{wolfram} of the partition function is given by
 \[
 p(m)\approx \frac{a\,e^{b\sqrt{m}}}{m}, \quad\mbox{where}\quad a=\frac{1}{4\sqrt{3}},\quad b=\frac{\pi\sqrt{2}}{\sqrt{3}} .
 \]
 Put $f(t)=e^{b\sqrt{t+1}}/(t+1)$ and $F(t)=\int_0^t f(x)\,dx$, yielding the convolution
 \[
 \int_0^t (t-x)f(x)\,dx=\int_0^t F(x)\,dx\approx \frac{2}{b}\int_0^t \frac{e^{b\sqrt{x+1}}}{\sqrt{x+1}}\,du\approx  \frac{4}{b^2} e^{b\sqrt{t}}.
 \]
 With $t\approx {n/2}$ we adjust the constants,  omitting the lower order term $p(m)$ when $n$ is even.
 \QED
\vv
The relatively low count warns us about possible inefficiency of the product procedure. However, this count applies only to canonical $P$ but they do not commute with permutations, and these would bring the factorial multiplier. Therefore, in spite of the initial simplicity the products of such basic matrices quickly become rich and complicated, which suggests that there might be a multiplication algorithm leading to an arbitrary D-orthogonal matrix.
 \begin{theorem}
 Every D-orthogonal matrix of dimension $n$ is indistinguishable from a product of $n$  D-orthogonal matrices of the form $I+\su\otimes\su$.
\end{theorem}
\Proof
It suffices to show that a D-orthogonal matrix $P$ admits a sequence $K_i=I+\su_i\otimes\su_i$ and a permutation matrix $\pi$  such that $K_n\cdots K_1 P=\pi$. To this end, note the implication
\[
\la\sp\sce\ra\cong 0\,\,\mbox{and}\,\, \la \sp\ra\cong 1\quad \imp \quad \Big(I+(\sp+\sce)\otimes (\sp+\sce)\Big)\,\sp\cong\sce.
\]
Whence $(I+\su\otimes \su)\,\sp\cong \si_k$ for an odd $\sp\neq \si$ and $p_k=0$, denoting $\su=\sp+\si_k$.\vv

Let us focus on the first column $\sp_1\neq \si$ of $P$ and create  an even vector $\su_1$ by the above recipe. Then the first column of $(I+\su_1\otimes \su_1)P$ contains a solitary $1$. Since the product is D-orthogonal hence this 1 is solitary in the corresponding row. The new columns are still odd and have even intersections with other columns, and there is no $\si$. Then, the procedure, repeated for the second (new) column, then for the third one, etc.,  without affecting the previously obtained single 1s, entails a permutation matrix $\pi$ at the end.
\QED

\section{Types of signatures}\label{sect:signs}
We will classify finite or infinite AC generators with respect to the distribution of signatures. The presence of one or more additional negative elements that commute with all others provides new means of controlling signatures. We will discuss this case in Section \ref{sign class}.
Recall that  $s_+$ or $s_-$ denote the counts of powers $\wek e^{\wek p}$ that are negative depending  upon the generator $\wek e$ being positive or negative, respectively. We will evaluate the counts.

\begin{proposition}\label{countsel}
Consider a pure anticommutative generator $\se$ of length $n$.
\begin{enumerate}
\item
We have $s_-=s_+$ if and only if $n$ is divisible by 4. In this case positive and negative generators are replaceable.
\item
The counts  are as follows:
\[
\begin{array}{rl}
\vspandex s_+&=2^{n-2} -2^{n/2-1} \left(\cos\Frac{n\pi }{4}+\sin\Frac{n\pi }{4}\right)\\
s_-&=2^{n-2} -2^{n/2-1} \left(\cos\Frac{n\pi }{4}-\sin\Frac{n\pi }{4}\right)\\
\end{array}
\]
\end{enumerate}
\end{proposition}
\Proof
We will use $b_q$ from \refrm{ss}. The 1834 Ramus identity (cf. \cite{Ram,Kon}) involved a period $p$ and the $p$\tss{th} root $\omega=e^{2\pi\imath/p}$ of 1:
\[
b_q(p)=\sum_j {n\choose q+pj}= \frac{1}{p} \sum_{k=1}^p \omega^{-qk}(1+\omega^k)^n.
\]
It follows by the binomial formula used in the right hand side and the fact that for $z=\omega^{r}$
\[
\sum_{k=1}^p z^k \neq 0 \quad\iff\quad p|r,
\]
in which case the sum equals $k$. However, we just need its very elementary version, because the period $p=4$ entails $\omega=\imath$:
\be\label{bq}
b_q=\frac{1}{4} \sum_{k=1}^4 \imath^{-qk}(1+\imath^k)^n= 2^{n-2}+\frac{c_q(n)}{4},
\quad\mbox{where}\quad c_q(n)=\frac{(1+\imath)^n}{\imath^q}+\frac{(1-\imath)^n}{\imath^{3q}}.
\ee
Hence
\be\label{cc}
c_1=\frac{(1+\imath)^n-(1-\imath)^n}{\imath}=-c_3,\quad
c_2=-(1+\imath)^n-(1-\imath)^n.
\ee
Therefore, $s_+=s_-$ if and only if $c_1=0$, i.e., $(1+\imath)^n=(1-\imath)^n$, or, in other words
\[
{\imath^n=\left(\frac{1+\imath}{1-\imath}\right)^n}=1\quad\iff\quad 4|n.
\]
So, consider $n=4m$. For $m=1$ let us choose  $ A=\set{1110,\,1101,\,1011,\,0111}$. Then the four triple products $\se^\sp$, $\sp\in A$,  anticommute, since $pq-\la \sp\sq\ra=9-2=7$ for $\sp,\sq\in A$. If $\se$ is negative then $\se^\sp$ are positive  since their signatures are $(-1)^{p(p+1)/2}=(-1)^6$.  If $\se$ is positive then $\se^\sp$ are negative  since their signatures are $(-1)^{p(p-1)/2}=(-1)^3$.
\vv

If $m>1$, we partition $\se$ into $m$ disjoint quadruples and apply the above signature changing process within each quadruple to bring all elements to the same signature. The new elements from disjoint quadruples anticommute because their binary marks $\sp$ and $\sq$ produce the form $pq-\la\sp\sq\ra=9-0$. This completes the proof of the first statement.\vv

Now, we combine the sums in \ref{ss} with \refrm{bq} and \refrm{cc}:
\[
c_1=2^{n/2-1}\sin\Frac{n\pi}{4}=-c_3,\quad c_2=-2^{n/2-1}\cos\Frac{n\pi}{4}.
\]
In the periodic trigonometric sequences $t_{\pm}=\cos\Frac{n\pi }{4}\pm\sin\Frac{n\pi }{4}$ of period 8 the radicals (or 0) appear only at some odd indices, so they simplify due to the factor $2^{n/2-1}$:
\[
\begin{array}{rccccccccccl}
n=                      &... &0,&1,&2,&3,&4,&5,&6,&7,&...\\
t_+: &...&1,&  \sqrt{2},& 1, &0,& -1,& -\sqrt{2},& -1,&0,&...\\
t_-:&...& 1,& 0,& -1,&-\sqrt{2},&-1,&  0,& 1, &\sqrt{2},&...\\
\end{array}
\]
This concludes the proof.
\QED
\vspace{5pt}

\vv

The sign matrix $C$ of a finite group with an anticommutative generator is unique. However, the signatures may be arbitrary and may change under  AC replacements. Let us denote by $N_-=N_-(\se)$ and $N_+=N_+(\se)$ the number of negative and positive elements, respectively, in an AC generator $\se$. If its length is $n$, so $N_-+N_+=n$.  Clearly, if $N_-(\se)=N_-(\se')$ then the generated groups are isomorphic. Yet, $N_-(\se)$ is not an isomorphism invariant and it is not immediately clear how the quantity behaves under {\bf admissible} replacement, $\se\mapsto \se^P$, where $P$ is DI and AC-preserving.
\begin{example}\rm Let us write $(m,k)=(N_+,N_-)$.

\begin{enumerate}
\item
The replacement $F\mapsto F(f_0)$ from Example \ref{addcolrow}.1, where an element $f_0$ is positive, changes the signatures as follows
\[
(m,k) \quad\leftrightarrow\quad (k+1,m-1).
\]
\item Let $n=4$ and  $\se=(e_1,e_2,e_3,e_4)$ be pure. Then $P=\Ic\!_4$ defines the admissible  replacement $\se'=\se^P=(e_2e_3e_4,\,e_1e_3e_4,\,e_1e_2e_4,\,e_1e_2e_3)$ with the signatures swapped:
\[
(4,0) \quad\leftrightarrow\quad (0,4).
\]
For the other configurations: $(3,1)\quad\leftrightarrow\quad (2,2)$, which is also covered by the previous replacement, $(1,3)\quad\leftrightarrow\quad (1,3)$, with $(1,3)\quad\leftrightarrow\quad (4,0)$. Thus we obtain two non-isomorphic signed groups, with equivalent distribution of $(+,-)$ elements: $\set{(4,0), (0,4), (1,3)}$ and $\set{(3,1), (2,2)}$.
 We will prove it in the general case.
\end{enumerate}
\end{example}
\begin{theorem}\label{signatures}
Let $n\ge 2$ be the length of an AC generator of a signed group of order $2^{n+1}$.  Up to an isomorphism,
\begin{enumerate}
\item if $n$ is even then there are two non-isomorphic groups;
\item if $n$ is odd then there are three non-isomorphic groups;
\item if $n=\infty$ then there is only one signed group.
\end{enumerate}
\end{theorem}
\Proof
There are $n+1$ possible replacements. We record the replacement within same group based on transformations in the above Example, either $m=N_+$ belongs to an arithmetic sequence of step 4, or being subject to the transformation $(m,k)\leftrightarrow (k+1,m-1)$.\vv

If $n=2$, two  groups $(2,0)\leftrightarrow (1,1)$ and $(0,2)$ are not isomorphic because pure generators yield different counts of negative elements by Proposition \ref{countsel}. So let $n\ge 3$.\vv

For $n=4k+2$ we obtain exactly two non-isomorphic groups. The first group has a pure negative generator while the second one has a pure positive generator, hence the total counts $N(\ca G)$ differ by Proposition \ref{countsel}. The negative count is isomorphism invariant.\vv
\[
\begin{array}{ll}
\begin{array}{cl}
(4i, 4(k-i)+2)      &\rule{17pt}{0pt}  i=0,...,k \\
\rule{4pt}{0pt}\updownarrow&\\
(4(k-i)+3,4i-1),    &\rule{17pt}{0pt} i=1,...,k
\end{array}&
\left.\vspandex\right\}\,\mbox{($2k+1$ replacements)}
\\
&\\
\begin{array}{cl}
(4i+2, 4(k-i))      &i=0,...,k \\
\updownarrow\rule{4pt}{0pt}&\\
 (4(k-i)+1,4i+1), &  i=0,...,k
\end{array}&
\left.\vspandex\right\}\,\mbox{($2k+2$ replacements)}
\end{array}
\]
In other words, the two types in the case of $n=2$ (mod 4) follow again the remainder of $N_+/4$:

\begin{center}
\begin{tabular}{ll}
type $\ca R_{n:2}(0,3)$: &if $N_+=0$ (mod 4) or $N_+=3$ (mod 4);\\
type $\ca R_{n:2}(1,2)$: &if $N_+=1$ (mod 4) or $N_+=2$ (mod 4).
\end{tabular}
\end{center}
In particular, the split $(N_+,N_-)=(0,n)$ characterizes $\ca R_{n:2}(0,3)$.\vv

For $n=4k+3$ we obtain at most three non-isomorphic groups:
\[
\begin{array}{ll}
\begin{array}{cl}
(4i, 4(k-i)+3)      &\rule{17pt}{0pt}  i=0,...,k \\
\rule{4pt}{0pt}\updownarrow&\\
\mbox{the same}   &
\end{array}&
\left.\vspandex\right\}\,\mbox{($k+1$ replacements)}\\
&\\
\begin{array}{cl}
(4i+1, 4(k-i)+2)      &i=0,...,k \\
\updownarrow\rule{4pt}{0pt}&\\
 (4(k-i)+3,4i), &  i=0,...,k
\end{array}&
\left.\vspandex\right\}\,\mbox{($2k+2$ replacements)}\\
&\\
\begin{array}{cl}
(4i+2, 4(k-i)+1)      &i=0,...,k \\
\updownarrow\rule{4pt}{0pt}&\\
 \mbox{the same}, &
\end{array}&
\left.\vspandex\right\}\,\mbox{($k+1$ replacements)}\\
\end{array}
\]
In other words,  for $n=3$ (mod 4) three types emerge
\begin{center}
\begin{tabular}{ll}
type $\ca R_{n:3}(0)$:    &if $N_+=0$ (mod 4);\\
type $\ca R_{n:3}(1,3)$: &if $N_+=1$ (mod 4) or $N_+=3$ (mod 4);\\
type $\ca R_{n:3}(2)$:   &if $N_+=2$ (mod 4);
\end{tabular}
\end{center}
They are not isomorphic. Indeed, by Proposition \ref{countsel} groups of the first and of the second type are not isomorphic because they have pure generators of opposite signatures. Also, each of them contains the subgroup of type $\ca R_{n-1:2}(0,3)$, i.e., one with generator $(N_+,N_-)=(0,n-1)$ of length $n-1$, which does not show among subgroups of a group of type $\ca R_{n:3}(2)$, because the maximum $N_-=4k+1=n-2$.

For $n=4k$  two groups emerge:
\[
\begin{array}{ll}
\begin{array}{cl}
(4i, 4(k-i))      &\rule{17pt}{0pt}  i=0,...,k \\
\rule{4pt}{0pt}\updownarrow&\\
(4(k-i)+1,4i-1),    &\rule{17pt}{0pt} i=1,...,k
\end{array}&
\left.\vspandex\right\}\,\mbox{($2k+1$ replacements)}
\\
&\\
\begin{array}{cl}
(4i+2, 4(k-i)-2)      &i=0,...,k-1 \\
\updownarrow\rule{4pt}{0pt}&\\
 (4(k-i-1)+3,4i+1), &  i=0,...,k-1
\end{array}&
\left.\vspandex\right\}\,\mbox{($2k$ replacements)}
\end{array}
\]
In other words, when $n=0$ (mod 4):

\begin{center}
\begin{tabular}{ll}
type $\ca R_{n:0}(0,1)$: &if $N_+=0$ (mod 4) or $N_+=1$ (mod 4);\\
type $\ca R_{n:0}(2,3)$: &if $N_+=2$ (mod 4) or $N_+=3$ (mod 4).
\end{tabular}
\end{center}

The lack of isomorphism follows again by a subgroup argument. That is, a group of type $\ca R_{n:0}(0,1)$ contains a subgroup $\ca R_{n-1:3}(0)$, represented by $(N_+,N_-)=(0,n-1)$, which is absent among subgroups of $\ca R_{n-1:3}(2,3)$ for which the maximum $N_-=4k-2=n-2$.\vv

The last case $n=4k+1$ follows similarly with three groups at hand.
\[
\begin{array}{ll}
\begin{array}{cl}
(4i, 4(k-i)+1)      &\rule{17pt}{0pt}  i=0,...,k \\
\rule{4pt}{0pt}\updownarrow&\\
(4(k-i)+2, 4i-1)   &\rule{17pt}{0pt}  i=1,...,k
\end{array}&
\left.\vspandex\right\}\,\mbox{($2k+1$ replacements)}\\
&\\
\begin{array}{cl}
(4i+1, 4(k-i))      &i=0,...,k \\
\updownarrow\rule{4pt}{0pt}&\\
\mbox{the same} &
\end{array}&
\left.\vspandex\right\}\,\mbox{($k+1$ replacements)}\\
&\\
\begin{array}{cl}
(4i+3, 4(k-i-1)+2)      &i=0,...,k-1 \\
\updownarrow\rule{4pt}{0pt}&\\
 \mbox{the same}, &
\end{array}&
\left.\vspandex\right\}\,\mbox{($k$ replacements)}\\
\end{array}
\]
Like before,  $n=1$ (mod 4) yields three types of groups:
\begin{center}
\begin{tabular}{ll}
type $\ca R_{n:1}(0,2)$:    &if $N_+=0$ (mod 4) or $N_+=2$ (mod 4);\\
type $\ca R_{n:1}(1)$: &if $N_+=1$ (mod 4);\\
type $\ca R_{n:1}(3)$:   &if $N_+=3$ (mod 4);
\end{tabular}
\end{center}
$\ca R_{n:1}(0,2)$ with a positive generator and $\ca R_{n:1}(1)$ with a negative generator are not isomorphic. Also,  each contains $\ca R_{n-1:0}(0,1)$, displaying a generator with $N-_=n-1$, as a subgroup in contrast to $\ca R_{n:1}(3)$ where the maximum $N_-=4k-2=n-3$.
\vv

Let us now turn to the proof of the third statement.
Using again the classification $(N_+,N_-)$, we first reduce all possible groups to four cases: \vv

\hspace{60pt} 1. $(0,\infty), (4,\infty),...\quad\leftrightarrow\quad (\infty, 3),\,(\infty, 7),...$;

\hspace{60pt} 2. $(1,\infty), (5,\infty),...\quad\leftrightarrow\quad (\infty, 0),\,(\infty, 4),...$;

\hspace{60pt} 3. $(2,\infty), (6,\infty),...\quad\leftrightarrow\quad (\infty, 1),\,(\infty, 5),...$;

\hspace{60pt} 4. $(3,\infty), (7,\infty),...\quad\leftrightarrow\quad (\infty, 2),\,(\infty, 6),...$,
\vv
\noindent each represented by a sequence that begins with 0, 1, 2, or 3 positive elements. Then, focusing on the first four elements, we use the type $\ca R_{4:0}(0,1)$ to merge the cases 1 and 2 into one case, represented by $(0,\infty)$. Then we use $\ca R_{4:0}(2,3)$  to merge the cases 3 and 4 to one case, represented by $(2,\infty)$.\vv

Then we cut off the first five elements and consider $(0,5)$ versus $(2,3)$ to see that $\ca R_{5:1}(0,2)$ makes one group, while the rest of the sequence, which is negative, remains unaffected.
\QED

\section{Classification of finite signed systems}\label{sect:parts}
 We consider a signed group $\ca G$ together with it basic subsets $E$, possibly ordered $E\leftrightarrow \se$.   It is convenient to refer to the  0-1 equivalent $D$ of the AC-matrix $C$, cf.\ \refrm{CD}. For a time being we disregard the signatures of elements because they do not affect the AC-matrix. By convention, we consider a singleton $\set{g}$ and $\emptyset$ as both commutative and anticommutative.
 We denote by $c^-=c^-(E)$ (also called  the `{\bf AC-count}') the count of the negative signs in the AC-matrix $C$ (or 1s in $D$) of the group generated by the generator $E$. The quantity is a group invariant, thus different counts yield non-isomorphic groups. We will see that the inverse implication is conditionally true, i.e., equal counts imply an isomorphism provided that the signature patterns agree.\vv

 A commutative signed generator happens if and only if the generated group is commutative. A more interesting situation occurs when some elements anticommute. Although arbitrary or even random signs may be assigned  to a generator, possibly entailing chaos of signs  in the AC-matrix of the generated group, yet the opposite happens. We will see that  quite `orderly' replacements exist that enjoy clear sign patterns. In particular, a subset of elements commuting with all others can be put aside.

 \begin{example}\rm
With $D$ of size $n\times n$, if 1s are placed at random  in the upper triangle $D\triu$ then the probability of obtaining at least one element commuting with all others is
\[
1-\prod_{k=1}^{n-1}\left(1-2^{-k}\right)\to 1-\phi(1/2)\approx 0.7112\mbox{ as $n\to\infty$}.
\]
Indeed, we may view $D$ as the adjacency matrix of a graph with $n$ vertices. If $p_n$ denotes the probability that the graph is connected, then $p_2=2^{-1}=1-2^{-1}$ and, conditioning on the first vertex, $p_n=\left(1-2^{n-1}\right) p_{n-1}$ for $n\ge 3$.\QED
 \end{example}
\vv

 \subsection{Equivalents of AC generators}
 Below we will give meaning to two actions: {\bf \em integration} of smaller (desirably simple) matrices of the given type to form a larger matrix of the same type, and  the inverse action of {\bf \em disintegration}.
  \begin{proposition}\label{example}
{~}
  \begin{enumerate}
 \item When $n=2$ or $n=3$ then at least one pair of anticommuting elements yields an anticommutative generator.
 \enuproof
     {The case $n=2$ is trivial. Let $n=3$ and let $D\triu$ have either a single 0 or a single 1, e.g., neglecting equivalent permutations:
 \[
 \bmat \cdot & 0 & 1\\ \cdot & \cdot & 1 \\ \cdot & \cdot & \cdot \emat \quad\mapsto\quad
 (e_1,\,e_2e_3,e_3),\quad
 \bmat \cdot & 1 & 0\\ \cdot & \cdot & 0 \\ \cdot & \cdot & \cdot \emat \quad\mapsto\quad
 (e_1,\,e_2,\,e_1 e_2 e_3).
 \]}
 \item
A chain-like sequence $\su$ such that $u_{k}\circ u_{k+1}=-1$ and  $u_i\circ u_j=1$ when $|i-j|\ge 2$ can be replaced by the AC-generator $e_k=u_1\cdots u_k$. Conversely, $u_k= e_{k-1} e_k$  (here $e_0=1$).
\enuproof
{Indeed, the needed properties follow by inspection.}
\item\label{doubletons}
An AC generator $E$ of even or infinite size and a chain (i.e., the union) of commuting AC doubletons $D_k=\set{d_{2k-1},d_{2k}}$ (i.e., $D_i\circ D_j=1$ for $i\neq j$ and $d_{2k-1}\circ d_{2k}=-1$) are mutually replaceable.
\enuproof
{Indeed, the repetitive replacement $E'=D\cup DE$, where $D=\set{e_1,e_2}$, yields a chain. Explicitly, the products $f_k=e_1\cdots e_k$ entail $E'=\bigcup_n D_n$ with  $D_n=\set{d_{2n-1},d_{2n}}$, where
\[
d_1=e_1,\,d_2=e_2, \quad \mbox{and for $n\ge 2$,}\quad
  d_{2n-1}=f_{2n-2}\,e_{2n-1},\quad d_{2n}=f_{2n-2}\,e_{2n}.
\]
The replacement is self-invertible:
\[
e_{2n-1}=d_1\cdots d_{2n-2}\,d_{2n-1},\,e_{2n}=d_1\cdots d_{2n-2}\,d_{2n}.
\]}
\item\label{arb}
In particular, a finite even or infinite  AC generator $E$ can be replaced by a generator $E'$ that admits a partition $E'=F_1\cup F_2\cup\cdots$ with mutually commuting AC components of finite even or infinite size. For an infinite $E$ the number of components can be either finite or infinite.
\enuproof
{Indeed, the aforementioned {\bf\em integrations} of doubletons or {\bf\em disintegrations} into doubletons can be combined at will to obtain an arbitrary described replacement.
 }
\item  Let $E=K\cup M$, where $K\circ M=1$, $K$ is finite and AC, and $M$ is commutative. If the size of $K$ is even and $M\neq\emptyset$, then $K$  can be increased by one, while if $K$ is of odd size it may be decreased by one.
    \enuproof
{Indeed, for $K=\set{k_1,\dots, k_{2j}}$ and $m_0\in M$, we define $k_0=k_1\cdots k_{2j}m_0$. Since $k_0\circ K=-1$ then we can replace $K'=K\cup \set{k_0}$ and $M'=M\setminus\set{m_0}$. If $K=\set{k_0,k_1,\dots, k_{2j}}$, then we augment $M$ by the total product $k_0k_1\cdots k_{2j}$, which commutes with all of $K$, and remove $k_0$ from $K$.\QED
}\end{enumerate}
 \end{proposition}

An AC generator is shown in the diagram below. Possible additional elements, independent of and commuting with the generator, are not displayed.

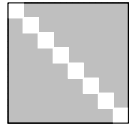
\begin{center}
\begin{tikzpicture}[scale=0.2]
\draw[help lines,dotted] (0,0) grid (8,8);
\draw [fill=lightgray] (0,0) rectangle (8,8);
\draw [fill=white,white] (0,7) rectangle (1,8);
\draw [fill=white,white] (1,6) rectangle (2,7);
\draw [fill=white,white] (2,5) rectangle (3,6);
\draw [fill=white,white] (3,4) rectangle (4,5);
\draw [fill=white,white] (4,3) rectangle (5,4);
\draw [fill=white,white] (5,2) rectangle (6,3);
\draw [fill=white,white] (6,1) rectangle (7,2);
\draw [fill=white,white] (7,0) rectangle (8,1);

\draw[black,thin] (0,8) --(8,8) -- (8,0)  -- (0,0) -- (0,8);
\end{tikzpicture}
\end{center}

An AC generator is replaceable as described in Points 2 through 4:
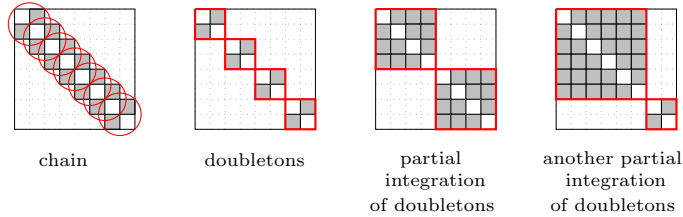
\begin{center}
\begin{tikzpicture}[scale=0.2]
\draw[help lines,dotted] (0,0) grid (8,8);
\draw [fill=lightgray] (1,7) rectangle (2,8);
\draw [fill=lightgray] (2,6) rectangle (3,7);
\draw [fill=lightgray] (3,5) rectangle (4,6);
\draw [fill=lightgray] (4,4) rectangle (5,5);
\draw [fill=lightgray] (5,3) rectangle (6,4);
\draw [fill=lightgray] (6,2) rectangle (7,3);
\draw [fill=lightgray] (7,1) rectangle (8,2);
\draw [fill=lightgray] (0,6) rectangle (1,7);
\draw [fill=lightgray] (1,5) rectangle (2,6);
\draw [fill=lightgray] (2,4) rectangle (3,5);
\draw [fill=lightgray] (3,3) rectangle (4,4);
\draw [fill=lightgray] (4,2) rectangle (5,3);
\draw [fill=lightgray] (5,1) rectangle (6,2);
\draw [fill=lightgray] (6,0) rectangle (7,1);

\draw[black,ultra thin] (0,8) --(8,8) -- (8,0)  -- (0,0) -- (0,8);
\draw [red, thin] (1,7) circle [radius=1.4142];
\draw [red, thin] (2,6) circle [radius=1.4142];
\draw [red, thin] (3,5) circle [radius=1.4142];
\draw [red, thin] (4,4) circle [radius=1.4142];
\draw [red, thin] (5,3) circle [radius=1.4142];
\draw [red, thin] (6,2) circle [radius=1.4142];
\draw [red, thin] (7,1) circle [radius=1.4142];
\node[right] at (1,-2) {\scriptsize chain};
\draw[help lines,dotted] (12,0) grid (20,8);
\draw [fill=lightgray] (13,7) rectangle (14,8);
\draw [fill=lightgray] (15,5) rectangle (16,6);
\draw [fill=lightgray] (17,3) rectangle (18,4);
\draw [fill=lightgray] (19,1) rectangle (20,2);
\draw [fill=lightgray] (12,6) rectangle (13,7);
\draw [fill=lightgray] (14,4) rectangle (15,5);
\draw [fill=lightgray] (16,2) rectangle (17,3);
\draw [fill=lightgray] (18,0) rectangle (19,1);
\node[right] at (12,-2) {\scriptsize doubletons};
\draw[black,ultra thin] (12,8) --(20,8) -- (20,0)  -- (12,0) -- (12,8);
\draw[red,thick] (12,8) --(14,8) -- (14,4) -- (18,4) -- (18,0) -- (20,0) -- (20,2) -- (16,2) -- (16,6) -- (12,6) -- (12,8);
\draw[help lines,dotted] (24,0) grid (32,8);
\draw [fill=lightgray] (25,7) rectangle (26,8);
\draw [fill=lightgray] (26,7) rectangle (27,8);
\draw [fill=lightgray] (27,7) rectangle (28,8);
\draw [fill=lightgray] (24,5) rectangle (25,6);
\draw [fill=lightgray] (25,5) rectangle (26,6);
\draw [fill=lightgray] (27,5) rectangle (28,6);
\draw [fill=lightgray] (29,3) rectangle (30,4);
\draw [fill=lightgray] (30,3) rectangle (31,4);
\draw [fill=lightgray] (31,3) rectangle (32,4);
\draw [fill=lightgray] (28,1) rectangle (29,2);
\draw [fill=lightgray] (29,1) rectangle (30,2);
\draw [fill=lightgray] (31,1) rectangle (32,2);
\draw [fill=lightgray] (24,6) rectangle (25,7);
\draw [fill=lightgray] (26,6) rectangle (27,7);
\draw [fill=lightgray] (27,6) rectangle (28,7);
\draw [fill=lightgray] (24,4) rectangle (25,5);
\draw [fill=lightgray] (25,4) rectangle (26,5);
\draw [fill=lightgray] (26,4) rectangle (27,5);
\draw [fill=lightgray] (28,2) rectangle (29,3);
\draw [fill=lightgray] (30,2) rectangle (31,3);
\draw [fill=lightgray] (31,2) rectangle (32,3);
\draw [fill=lightgray] (28,0) rectangle (29,1);
\draw [fill=lightgray] (29,0) rectangle (30,1);
\draw [fill=lightgray] (30,0) rectangle (31,1);
\draw[black,ultra thin] (24,8) --(32,8) -- (32,0)  -- (24,0) -- (24,8);
\draw[red,thick] (24,8) --(28,8) -- (28,0) -- (32,0) -- (32,4) -- (24,4) -- (24,8);
\node[right] at (25,-2) {\scriptsize partial};
\node[right] at (24,-3.5) {\scriptsize integration};
\node[right] at (23,-5) {\scriptsize of doubletons};
\draw[help lines,dotted] (36,0) grid (44,8);
\draw [fill=lightgray] (37,7) rectangle (38,8);
\draw [fill=lightgray] (38,7) rectangle (39,8);
\draw [fill=lightgray] (39,7) rectangle (40,8);
\draw [fill=lightgray] (40,7) rectangle (41,8);
\draw [fill=lightgray] (41,7) rectangle (42,8);

\draw [fill=lightgray] (36,6) rectangle (37,7);
\draw [fill=lightgray] (38,6) rectangle (39,7);
\draw [fill=lightgray] (39,6) rectangle (40,7);
\draw [fill=lightgray] (40,6) rectangle (41,7);
\draw [fill=lightgray] (41,6) rectangle (42,7);
\draw [fill=lightgray] (36,5) rectangle (37,6);
\draw [fill=lightgray] (37,5) rectangle (38,6);
\draw [fill=lightgray] (39,5) rectangle (40,6);
\draw [fill=lightgray] (40,5) rectangle (41,6);
\draw [fill=lightgray] (41,5) rectangle (42,6);
\draw [fill=lightgray] (36,4) rectangle (37,5);
\draw [fill=lightgray] (37,4) rectangle (38,5);
\draw [fill=lightgray] (38,4) rectangle (39,5);
\draw [fill=lightgray] (40,4) rectangle (41,5);
\draw [fill=lightgray] (41,4) rectangle (42,5);
\draw [fill=lightgray] (36,3) rectangle (37,4);
\draw [fill=lightgray] (37,3) rectangle (38,4);
\draw [fill=lightgray] (38,3) rectangle (39,4);
\draw [fill=lightgray] (39,3) rectangle (40,4);
\draw [fill=lightgray] (41,3) rectangle (42,4);
\draw [fill=lightgray] (36,2) rectangle (37,3);
\draw [fill=lightgray] (37,2) rectangle (38,3);
\draw [fill=lightgray] (38,2) rectangle (39,3);
\draw [fill=lightgray] (39,2) rectangle (40,3);
\draw [fill=lightgray] (40,2) rectangle (41,3);
\draw [fill=lightgray] (43,1) rectangle (44,2);
\draw [fill=lightgray] (42,0) rectangle (43,1);
\draw[black,ultra thin] (36,8) --(44,8) -- (44,0)  -- (36,0) -- (36,8);
\draw[red,thick] (36,8) --(42,8) -- (42,0) -- (44,0) -- (44,2) -- (36,2) -- (36,8);
\node[right] at (34.5,-2) {\scriptsize another partial};
\node[right] at (36,-3.5) {\scriptsize integration};
\node[right] at (35,-5) {\scriptsize of doubletons};
\end{tikzpicture}
\end{center}

In Point 5 we toggle between odd and even anticommutative portions of generators:

\begin{center}
\begin{tikzpicture}[scale=0.2]
\draw[help lines,dotted] (0,0) grid (8,8);
\draw [fill=lightgray] (0,3) rectangle (5,8);
\draw [fill=white,white] (0,7) rectangle (1,8);
\draw [fill=white,white] (1,6) rectangle (2,7);
\draw [fill=white,white] (2,5) rectangle (3,6);
\draw [fill=white,white] (3,4) rectangle (4,5);
\draw [fill=white,white] (4,3) rectangle (5,4);

\draw[black,thin] (0,3) --(5,3) -- (5,8)  -- (0,8) -- (0,3);
\draw[black,thin] (5,0) --(8,0) -- (8,3)  -- (5,3) -- (5,0);
\draw[black,thin] (0,8) --(8,8) -- (8,0)  -- (0,0) -- (0,8);
\node[left] at (-.5,5.5) {\scriptsize odd $K$};
\node[left] at (-.5,1.5) {\scriptsize $M$};
\node[right] at (8.5,4) {\scriptsize $\leftrightarrow$};
\draw[help lines,dotted] (12,0) grid (20,8);
\draw [fill=lightgray] (12,4) rectangle (16,8);
\draw [fill=white,white] (12,7) rectangle (13,8);
\draw [fill=white,white] (13,6) rectangle (14,7);
\draw [fill=white,white] (14,5) rectangle (15,6);
\draw [fill=white,white] (15,4) rectangle (16,5);

\draw[black,dashed,thin] (12,3) --(17,3) -- (17,8)  -- (12,8) -- (12,3);
\draw[black,thin] (12,4) --(16,4) -- (16,8)  -- (12,8) -- (12,4);
\draw[black,thin] (16,0) --(20,0) -- (20,4)  -- (16,4) -- (16,0);
\draw[black,thin] (12,8) --(20,8) -- (20,0)  -- (12,0) -- (12,8);
\node[right] at (20.5,5.5) {\scriptsize even $K'$};
\node[right] at (20.5,1.5) {\scriptsize $M'$};
\end{tikzpicture}
\end{center}

 \begin{remark}\rm{~}
 \begin{enumerate}
 \item
 The aforementioned replacements can be expressed in terms of matrix transformations $\se\mapsto {^P\se}$ (see the paragraph above \refrm{V}). For example, the replacement in \refrm{Ff} corresponds to $P=C_1$ from  Example \ref{addcolrow}.1.
 \item
 Commuting positive doubletons appeared before, isomorphically defined by \refrm{action}. Negative doubletons required `complexification' as shown below \refrm{action}.
 \QED
 \end{enumerate}
\end{remark}

\subsection{Orderly partitions}
We assume that all sequences or corresponding sets appearing below are basic with at least two anticommuting elements.
\begin{proposition}\label{CvsAC}
Let an element $g$ be independent of a finite anticommutative basic set $F$. Then there are three distinct possibilities:
\begin{enumerate}
\item
$g$ commutes with $F$,
\item
$g$ anticommutes with at least one element of $F$, in which case  $F\cup\set{g}$ can be replaced
\begin{enumerate}
\item either by an anticommutative generator,
\item or by $F'\cup\set{g'}$, where $F'$ is anticommutative and $g'$ commutes with $F'$.
\end{enumerate}
\end{enumerate}

\end{proposition}

\Proof
 The element $g$ entails a partition of $F=F_a\cup F_c$, where $g\circ F_a=-1$  and $g\circ F_c=1$. Recall (Example \ref{addcolrow}.1) that for $f_0\in F$
 the replacement
$F(f_0)=\set{f_0}\cup \set{f_0f:f\in F,\,f\neq f_0}$ preserves anticommutativity. In addition we observe that
\be\label{F2}
\mbox{\em If $F_a=\set{f_0}$ then $g$ anticommutes with the anticommutative generator $F(f_0)$}.
\ee
Suppose that $g$ anticommutes with at least one element of $F$, i.e., $F_a\neq\emptyset$. If $F_a$ is a singleton, then  \refrm{F2} yields case (a). Suppose there are at least two distinct elements $f_0,\,f_1\in F_a$. Replace $F$ by $F'=F(f_0)$ and $g$ by $g'=f_0f_1g$. Then we check that $|F'_a|=|F_a|-2$, preserving parity:
 \be\label{reduction}
g\,\circ\,
\set{\underset{-}{f_0},\,\underset{-}{f_1},...\underset{-}{f_a}...\big| ...\underset{+}{f_c}...}
\quad\mapsto\quad
 f_0f_1g\,\circ\, \set{\underset{+}{f_0},\,\underset{+}{f_0f_1},...\underset{-}{f_0f_a}...\big|...\underset{+}{f_0f_c}...},
\ee
with values of the commutativity function marked beneath the elements. Therefore,  by recursion we can reduce $F_a$ either to the empty set, yielding (b), or to a singleton and then we use \refrm{F2} to arrive in case (a).
\QED

\begin{remark}\rm
A partial extension is valid for infinite $F$ when $F_a$ or $F_c$ is finite:
\begin{enumerate}
\item If  $F_a$ is finite then $F\cup\set{g}$ can be replaced
\begin{enumerate}
\item  either by $F'\cup\set{g'}$ with an AC set $F'$ and $g'\circ F'=1$, for even $|F_a|$,
\item
or by an AC set $F'$, for odd $|F_a|$.
\end{enumerate}
\item If $F_c$ is finite and $F_a\neq \emptyset$, then w.l.o.g. we may assume that $F_a$ is finite.
\end{enumerate}
The proof for a finite $F_a$ mimics the corresponding proof above. The case of finite $F_c$  follows by  switching to  $F(f)$, where $f\in F_a$.
\QED
\end{remark}

Recall the simple yet useful property \refrm{F2}.
The equivalence relation $\ca G(\se)=\ca G(\se')$ (i.e., sequences are mutually replaceable) between basic sequences extends to the relation of partial order $\ca G(\se)\subset  G(\se')$  between (equivalence classes of) basic sequences. The first `orderly pattern' will appear as follows:
\be\label{III}
E\mapsto E'\leftrightarrow D'=
\bmat
\Ic    &  O      & O      & \cdots & O     \\
O      & \Ic     & O      & \cdots & O     \\
\vdots & \vdots  & \ddots & \vdots & \vdots\\
O      & \cdots  & O      & \Ic    & O     \\
O      & \cdots  & O      & O      & O     \\
\emat
\ee
with $\Ic$'s of some varying sizes at least 2 along the diagonal. Suppose that all sizes are even. Then some or all of $\Ic$'s  may be turned to be odd if the commuting part appearing in the lower right corner is large enough, which can be always assumed by augmenting the system. See Remark \ref{evenodd} for more details.
\begin{theorem}\label{deco}
Let $E$ be a signed basic set of length $n\ge 2$ with at least one pair of anticommuting elements. Then for some $k=1,\dots,n$ there exists a replacement $E'$ of $E$ and its disjoint partition  $E'=F_0\cup F_1\cup\cdots \cup  F_k$ such that for each for $j\ge 1$ the set $F_j$  is anticommutative of at least length 2, $F_0$ is commutative, and $F_j\circ F_j=1$ for $j\neq j',\, 0\le j,\,j'\le k$.
\end{theorem}
\Proof
If a finite signed basic set $E$ has at least one anticommutative pair, then there exists a nonempty anticommutative basic  $F$ of maximal size. If $\ca G(F)=\ca G(E)$, we finish with $F_1=F$ and $F_0=\emptyset$. If $\ca G(F)\neq \ca G(E)$, we put $G=E\setminus F$ and then for any $g\not\in \ca G(F)$ we arrive in  Case 1 of Proposition \ref{CvsAC}, since $F$ is maximal. That is, $F\cup\set{g}$ is replaced by $F'\cup\set{g'}$ with an anticommutative $F'$ and $g'\circ F'=1$. If $G$ is commutative and commutes with $F$, then we are done with $F_1=F$ and $F_0=G$.\vv

Otherwise, consider an  anticommutative basic set $G_2\subset \ca G(E)\setminus \ca G(F)$, maximal in size. If $G_2\circ F=1$, we put $F_1=F$ and $F_2=G_2$.  If $G_2$ partially anticommutates with $F$, so each its element anticommutates with an even number of elements of $F$, since Proposition \ref{CvsAC} enforces Case 1 in view of the maximality of $F$. Then, browsing $G$, we skip elements $g$ that commute with $F$ (i.e., that `even number' is 0), and attend its subset  $G'$ of elements $g'$ admitting anticommutants in $F$. Then, for each of these elements $g'$ we apply \refrm{reduction} as many times as necessary to replace $F\cup\set{g'}$ by $F'\cup\set{g'}$. Note that the former elements $g$ still satisfy $g\circ F'=1$ and $g\circ g'=-1$. This part of the algorithm ends when we reach the last element in $G$, finishing with $F_1=F'$ and $F_2=G\cup G'$.\vv

Suppose that we arrived in a partition $E=E_k\cup G$, where $E_k=F_1\cup\cdots \cup F_k$  with anticommmutative maximal sets $F_i$ that commute among themselves.  Suppose that the partition is not final. We need to show how to find a subset of $G$ to create $F_{k+1}$. So, let $G_{k+1}$ be a maximal in size anticommutative set independent of $E_k$ that partially anticommutates with $E_k$. We skip its subset that fully commutes with $E_k$, and attend elements, one by one, that have aniticommutants in $E_k$. For each of these elements we apply procedure \refrm{reduction} with respect to consecutive $F_k$'s bearing anticommutants. We observe that the replacements $F\mapsto F(f)$ preserves the commutativity among $F_k$'s, and at the same time preserves the anticommutativity among $g$'s. This leads to $F_{k+1}$. \vv

The algorithm must end after finitely many steps with the sought-for partition.
\QED
\vv

Let us summarize our findings, allowing some redundancy, to exhibit two extremal patterns within the framework of \refrm{III}: on one extreme a single $\Ic$ and on the other doubletons yielding $2\times 2$ tiny matrices $\Ic={\left[\genfrac{}{}{0pt}{}{01}{10}\right]}$.

 \begin{theorem}\label{fullor2}
 Let $E$ be a basic signed set of size $n$ with at least two anticommuting elements. Then the following replacements exist.
 \begin{enumerate}
 \item
 There exists a replacement and its partition $E'=K\cup M$ such that $K$ is anticommutative, $M$ is commutative and may be empty, and $K$ and $M$ commute.
 \item
 Let $k$ denote the size of $K$ and $m$ denote the size of $M$. Then we may replace
     \[
     \begin{array}{ll}
     (k,m)\mapsto (k-1,m+1)&\mbox{if $k$ is odd}\\
     (k,m)\mapsto (k+1,m-1)&\mbox{if $k$ is even and $m\ge 1$}\\
     \end{array}
     \]
 \item There is a $K$ with the maximum even size $k=2j$. In this case there is a replacement and partition $K'=K_1\cup\cdots\cup K_l$ into the union of anticommutative doubletons that commute with each other and with $M$. Also, the AC-count is $c^-(E)=2^{2n-1}\left(2^{j}-1\right)$.
     \end{enumerate}
 \end{theorem}

\Proof
Consider any $F_i=\set{f_1,f_2,\dots}$ in the partition stated in Theorem \ref{deco}. The replacement of $F_i$:
\[
(f_1,\,f_2,...f,...)\mapsto (f_1,\,f_2,...,f_1f_2f,...)=(f_1,f_2)\,(...f_1f_2f,...)
\]
yields the partition of $F_i=F_{i1}\cup F_{i2}$ into two commuting anticommutative sets that still commute with all other sets $F_j$, $j\neq i$. We repeat the procedure while the size of $F_{i2}$ is at least 3, ending with a partition of $F_i$ into doubletons with desired properties and perhaps a single leftover that commutes with everything. In the latter case we augment $F_0$ by that singleton.\vv

The repetition of the algorithm for all components of the original partition entails Case 3, with a sequence of commuting doubletons. Now, given an anticommuting doubleton $\set{f_1,f_2}$ and an anticommuting $F'=\set{f'}$, and both commute, we replace their union by the anticommutating $\set{f_1,f_2,f_1f_2f': f'\in F'}$. The replacement does not change the commutativity with the remaining elements. Thus we end up with a $K=\set{f_1,\dots,f_k}$ of an even size $k$, which proves the first statement.
\vv

However, we may still modify $K$, as described in the second statement. If $k<n$, i.e., there is at least one element $f_0$ that commutes with all $f_k$, then $f_0f_1\cdots f_k$ anticommutes with $K$, which increases the size of $K$ by 1. The inverse replacement reduces the size of $K$ by 1.
\QED
\vspace{5pt}

 \begin{theorem}\label{quickcounts}
 Let $\se=(\se_1,\dots,\se_k)$ be a basic sequence of length $n$, consisting of subsequences $\wek e_i$ of length $\ell_i$ (so $\ell_1+\cdots+\ell_k=n$), such that $\se_i\circ\se_j=1$ for $i\neq j$. Let a commutative basic sequence $\sd$ of length $m$ commute with $\se$. Denote by $c_i=c^-(\se_i)$, the AC-count of the group generated by $\se_i$. Then
 \[
 c^-(\se)=\frac{1}{2}\left(2^{2n}-\prod_{i=1}^k \left(2^{2\ell_i}-2 c_i\right)\right),\qquad c^-(\se,\sd)=2^{2m}\,c^-(\se).
 \]
 In particular,
 \begin{enumerate}
 \item for an even $n$ and the entire AC $\wek e$, or equivalently, for  mutually commutative  pairs of AC anticommutating elements $\wek e_i$, $i=1,\dots, k$,
 $
 c^-(\se)=2^{n-1}\left(2^{n}-1\right)
 $
 \item for an odd $n$ and the entire  AC $\wek e$, the count equals
 $
 c^-(\se)=2^{n-2}\left(2^{n-1}-1\right).
 $
 \end{enumerate}
 \end{theorem}
 \Proof
 Consider two arrays $P=[\sp_1,\dots,\sp_k],\,Q=[\sq_1,\dots,\sq_k]$ and the sign
 \[
 s(P,Q)= \se_1^{\sp_1}\cdots \se_k^{\sp_k}\circ \se_1^{\sq_1}\cdots \se_k^{\sq_k}=
 (-1)^{C(\sp_1,\sq_1)+\cdots+C(\sp_k,\sq_k)}.
 \]
 Then AC-counts equal
 \[
c^-(\wek e_i)=\frac{1}{2}\sum_{\wek p_i,\wek q_i\in\BD_{\ell_i}}\left(1-s_i(\wek p_i,\wek q_i)\right),\quad \frac{1}{2}\,\sum_{P,Q}\left( 1-s(P,Q)  \right)=  \frac{1}{2}\,\sum_{P,Q}\left( 1-s_1\cdots s_k  \right),
 \]
 which yields the sought-for formula. Denote for the sake of brevity $N=c^-(\se)$. Then, for the augmented sequences $(\se,\sd)$ the AC-count follows from the diagram
  \vspace{5pt}

 \begin{center}
 \begin{tabular}{c||c|c|c}
  &$\se^{\sp}$  & $\sd^{\sq}$ &
  $\se^{\sp}\sd^{\sq}$, $q\neq 0$                               \\ \hline \hline
 $\se^{\sp}$           &     $N$       & $\cdot$ & $(2^m-1)\,N$ \\ \hline
 $\sd^{\sq}$           &   $\cdot$     & $\cdot$ &  $\cdot$     \\ \hline
  $\se^{\sp}\sd^{\sq}$, $q\neq 0$
                       &  $(2^m-1)\,N$ & $\cdot$ & $(2^m-1)^2\,N$
 \end{tabular}
 \end{center}
 \vspace{5pt}

For the pairs, the counts $s_i=6$ entail the special case when $n$ is even and $\wek e_i$ are commuting AC doubletons. \vv

Finally, we invoke Theorem \ref{fullor2} that allows us to pool doubletons together or make a replacement with a chain of doubletons, and an odd AC generator may be reduced by one, leaving a commuting element outside.
 \QED
\begin{remark}\label{evenodd}\rm
Recall partition \refrm{III}. While even commuting AC basic sequences can be integrated into one AC sequence, whose size is the sum of the sizes of parts, this does not occur when at least one of the sizes is odd. Each odd sequence must be first reduced by one to make it even, and only then the resulting even sequences can be integrated into one. For example, consider the sizes 9:7:2:(0), splitting a sequence of length 18, where the size of the commutative sub-generator commuting with all others appears in parentheses. Then we reduce two first sequences to arrive at the ratio 8:6:2:(2), setting aside two extra elements commuting with all. Next, we integrate the three even sequences into one AC sequence of length 16, i.e., 16:(2), which could be disintegrated at will, e.g., to 4:4:4:4:(2) or 12:4:(2), or else. Then the two extra elements can be added to some sequences, yielding, e.g., 5:5:4:4:(0) in the first case, or 13:5:(0) in the second case.\QED
\end{remark}
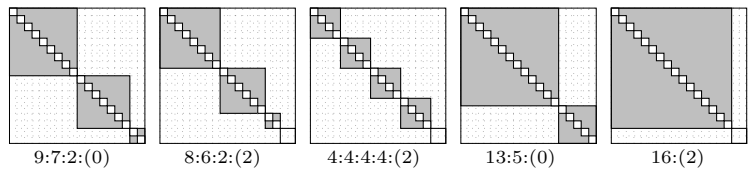
\begin{center}
\begin{tikzpicture}[scale=0.1]
\draw[help lines,dotted] (0,0) grid (18,18);
\draw [fill=lightgray] (0,9) rectangle (9,18);
\draw [fill=lightgray] (9,2) rectangle (16,9);
\draw [fill=lightgray] (16,0) rectangle (18,2);
\draw [fill=white] (0,17) rectangle (1,18);
\draw [fill=white] (1,16) rectangle (2,17);
\draw [fill=white] (2,15) rectangle (3,16);
\draw [fill=white] (3,14) rectangle (4,15);
\draw [fill=white] (4,13) rectangle (5,14);
\draw [fill=white] (5,12) rectangle (6,13);
\draw [fill=white] (6,11) rectangle (7,12);
\draw [fill=white] (7,10) rectangle (8,11);
\draw [fill=white] (8,9) rectangle (9,10);
\draw [fill=white] (9,8) rectangle (10,9);
\draw [fill=white] (10,7) rectangle (11,8);
\draw [fill=white] (11,6) rectangle (12,7);
\draw [fill=white] (12,5) rectangle (13,6);
\draw [fill=white] (13,4) rectangle (14,5);
\draw [fill=white] (14,3) rectangle (15,4);
\draw [fill=white] (15,2) rectangle (16,3);
\draw [fill=white] (16,1) rectangle (17,2);
\draw [fill=white] (17,0) rectangle (18,1);
\node[right] at (2,-2) {\scriptsize 9:7:2:(0)};

\draw[black,ultra thin] (0,18) --(18,18) -- (18,0)  -- (0,0) -- (0,18);
\draw[help lines,dotted] (20,0) grid (38,18);
\draw [fill=lightgray] (20,10) rectangle (28,18);
\draw [fill=lightgray] (28,4) rectangle (34,10);
\draw [fill=lightgray] (34,2) rectangle (36,4);
\draw [fill=white] (20,17) rectangle (21,18);
\draw [fill=white] (21,16) rectangle (22,17);
\draw [fill=white] (22,15) rectangle (23,16);
\draw [fill=white] (23,14) rectangle (24,15);
\draw [fill=white] (24,13) rectangle (25,14);
\draw [fill=white] (25,12) rectangle (26,13);
\draw [fill=white] (26,11) rectangle (27,12);
\draw [fill=white] (27,10) rectangle (28,11);
\draw [fill=white] (28,9) rectangle (29,10);
\draw [fill=white] (29,8) rectangle (30,9);
\draw [fill=white] (30,7) rectangle (31,8);
\draw [fill=white] (31,6) rectangle (32,7);
\draw [fill=white] (32,5) rectangle (33,6);
\draw [fill=white] (33,4) rectangle (34,5);
\draw [fill=white] (34,3) rectangle (35,4);
\draw [fill=white] (35,2) rectangle (36,3);
\node[right] at (22,-2) {\scriptsize 8:6:2:(2)};

\draw[black,ultra thin] (20,18) --(38,18) -- (38,0)  -- (20,0) -- (20,18);
\draw[black,thin] (36,0) --(38,0) -- (38,2)  -- (36,2) -- (36,0);
\draw[help lines,dotted] (40,0) grid (58,18);
\draw [fill=lightgray] (40,14) rectangle (44,18);
\draw [fill=lightgray] (44,10) rectangle (48,14);
\draw [fill=lightgray] (48,6) rectangle (52,10);
\draw [fill=lightgray] (52,2) rectangle (56,6);
\draw [fill=white] (40,17) rectangle (41,18);
\draw [fill=white] (41,16) rectangle (42,17);
\draw [fill=white] (42,15) rectangle (43,16);
\draw [fill=white] (43,14) rectangle (44,15);
\draw [fill=white] (44,13) rectangle (45,14);
\draw [fill=white] (45,12) rectangle (46,13);
\draw [fill=white] (46,11) rectangle (47,12);
\draw [fill=white] (47,10) rectangle (48,11);
\draw [fill=white] (48,9) rectangle (49,10);
\draw [fill=white] (49,8) rectangle (50,9);
\draw [fill=white] (50,7) rectangle (51,8);
\draw [fill=white] (51,6) rectangle (52,7);
\draw [fill=white] (52,5) rectangle (53,6);
\draw [fill=white] (53,4) rectangle (54,5);
\draw [fill=white] (54,3) rectangle (55,4);
\draw [fill=white] (55,2) rectangle (56,3);
\node[right] at (41,-2) {\scriptsize 4:4:4:4:(2)};
\draw[black,ultra thin] (40,18) --(58,18) -- (58,0)  -- (40,0) -- (40,18);
\draw[black,thin] (56,0) --(58,0) -- (58,2)  -- (56,2) -- (56,0);
\draw[help lines,dotted] (60,0) grid (78,18);
\draw [fill=lightgray] (60,5) rectangle (73,18);
\draw [fill=lightgray] (73,0) rectangle (78,5);
\draw [fill=white] (60,17) rectangle (61,18);
\draw [fill=white] (61,16) rectangle (62,17);
\draw [fill=white] (62,15) rectangle (63,16);
\draw [fill=white] (63,14) rectangle (64,15);
\draw [fill=white] (64,13) rectangle (65,14);
\draw [fill=white] (65,12) rectangle (66,13);
\draw [fill=white] (66,11) rectangle (67,12);
\draw [fill=white] (67,10) rectangle (68,11);
\draw [fill=white] (68,9) rectangle (69,10);
\draw [fill=white] (69,8) rectangle (70,9);
\draw [fill=white] (70,7) rectangle (71,8);
\draw [fill=white] (71,6) rectangle (72,7);
\draw [fill=white] (72,5) rectangle (73,6);
\draw [fill=white] (73,4) rectangle (74,5);
\draw [fill=white] (74,3) rectangle (75,4);
\draw [fill=white] (75,2) rectangle (76,3);
\draw [fill=white] (76,1) rectangle (77,2);
\draw [fill=white] (77,0) rectangle (78,1);
\node[right] at (62,-2) {\scriptsize 13:5:(0)};

\draw[black,ultra thin] (60,18) --(78,18) -- (78,0)  -- (60,0) -- (60,18);
\draw[help lines,dotted] (80,0) grid (98,18);
\draw [fill=lightgray] (80,2) rectangle (96,18);
\draw [fill=white] (80,17) rectangle (81,18);
\draw [fill=white] (81,16) rectangle (82,17);
\draw [fill=white] (82,15) rectangle (83,16);
\draw [fill=white] (83,14) rectangle (84,15);
\draw [fill=white] (84,13) rectangle (85,14);
\draw [fill=white] (85,12) rectangle (86,13);
\draw [fill=white] (86,11) rectangle (87,12);
\draw [fill=white] (87,10) rectangle (88,11);
\draw [fill=white] (88,9) rectangle (89,10);
\draw [fill=white] (89,8) rectangle (90,9);
\draw [fill=white] (90,7) rectangle (91,8);
\draw [fill=white] (91,6) rectangle (92,7);
\draw [fill=white] (92,5) rectangle (93,6);
\draw [fill=white] (93,4) rectangle (94,5);
\draw [fill=white] (94,3) rectangle (95,4);
\draw [fill=white] (95,2) rectangle (96,3);

\node[right] at (84,-2) {\scriptsize 16:(2)};

\draw[black,ultra thin] (80,18) --(98,18) -- (98,0)  -- (80,0) -- (80,18);
\draw[black,thin] (96,0) --(98,0) -- (98,2)  -- (96,2) -- (96,0);
\end{tikzpicture}
\end{center}
\begin{remark}\rm

 If $E$ is infinite, we encounter two ascending mutually commuting groups, with AC generators $K_n$ and commutative $M_n$.
Unfortunately, our finite algorithms, developed thus far, even augmented by a handful of infinite procedures, are not suitable for
infinite (countable) signed groups. In an infinite signed group $\ca G$  with infinitely many anticommuting elements we can find an ascending sequence of proper subgroups $\ca G_n$. The question is whether or not the union $\bigcup_n \ca G_n$ is proper. Should the answer be affirmative, its sub-generator $K$ would entail other patterns listed above or in Proposition \ref{example}.
The lack of  a quick answer is tied to our method of consecutive enlargements that are based on replacements, i.e., while groups grow their generators constantly change with no stabilization detected.\QED
\end{remark}

 \subsection{Partitions with signatures}\label{sign class}
An AC-doubleton always has a pure generator, either negative (of quaternion type), or positive (of Pauli type). Indeed, a mixed generator $(e_1,e_2)$ with $e_1^2=1,e_2^2=-1$ has the  positive replacement $(e_1,e_1e_2)$.  Invoke Theorem \ref{fullor2}. A commuting generator $M$ of size $m$ is either positive, or it has a replacement with any number $m'\ge 1$ of negative elements. Choose $m'=1$ and name that single positive or negative element $g$ (one stands for all) a `leftover'. \vv

 Consider a generator $E$ with an AC-pair and invoke the chain of doubletons described in Theorem \ref{fullor2}.3. Let $p$ and $q$ denote the number of positive and negative doubletons, respectively.
  That is, doubletons are either of Pauli type or of quaternion type. Assume that $M=\emptyset$ (so $E$ is even) or $M=\set{g}$ (so $E$ is odd). A larger positive commuting $M$ that commutes with everything else is irrelevant in our context, so we disregard it.\vv

  Therefore,  we can assign  the following temporary characteristic augmented by the signature $s=g^2$ of a leftover $g$ which appears only when the size $n$ of $E$ is odd:
 \def\PN#1#2{\Big\langle #1,#2\Big\rangle}
 \def\PNs#1#2#3{\Big\langle #1,#2;#3\Big\rangle}
\[
\begin{array}{ll}
\PN p q ,    &  \mbox{for even $n=2j$, $p+q=j$},\\
\PNs p q s,    &  \mbox{for odd $n=2j+1$, $p+q=j$}.
\end{array}
\]
We will see that there are factually either two or three characteristics, depending on parity.
We write $\PN p q=\PN {p'} {q'}$ if the corresponding generators are mutually replaceable. We will see below that only the parity of the number of negative (or positive) doubletons matters.  Thus, we arrive at the following taxa
\[
\PN 0 j,\, \PN 1{j-1}, \quad\mbox{or}\quad \PNs  0 j +,\, \PNs 1 {j-1} +,\,\PNs 0 j -,
\]
containing two groups for even $n=2j$, or  three groups for odd $n=2j+1$, respectively.

\begin{theorem}
If $n=2j$ is even, then $\PN p q=\PN {p'} {q'}$ iff $p\cong p'$. In other words, we obtain two types of groups:
\[
\PN 0 j=\PN 2 {j-2}=\cdots \quad \mbox{or}\quad \PN 1 {j-1} = \PN 3 {j-3}=\cdots
\]
i.e., the first group has an even number while the second group has an odd number of positive doubletons.
\vv
If $n=2j+1$ is odd, then for $s=-$ all decompositions $p+q=j$ are equivalent, and for $s=+$ two above groups appear.
\end{theorem}
\Proof
Let $n\ge 4$ be even. We examine the last step in the integration and disintegration $(e_1,e_2,...e...) \,\longleftrightarrow\, (e_1e_2)(...e_1e_2e...)$, described in the proof of Theorem \ref{fullor2} when a generator has length 4. As shown in Theorem \ref{signatures} there are
 two types $(4,0)\leftrightarrow (1,3)\leftrightarrow (0,4)$ and $(2,2)\leftrightarrow (3,1)$. We verify that for $n=4$
 \be\label{202}
 \PN 2 0\, \longleftrightarrow\, (2,2) \,\longleftrightarrow\, \PN 0 2.
 \ee
 The relation is just a symbolic expression of the explicit integration:
  \[
 (++)(++)\,\longleftrightarrow\,(++--) \quad\mbox{  and}\quad
 (--)(--)\,\longleftrightarrow \,(--++).
 \]
 Also, $ \PN 1 1\, \longleftrightarrow\, (4,0)$ since
 $(++)(--)\,\longleftrightarrow\, (++++)$.
 This means that $\PN 2 0= \PN 0 2$ and proves the statement for $n=4$. Let $n\ge 6$. Then relation \refrm{202} says that every pair of doubletons of the same parity can be replaced by a pair of doubletons of the opposite parity. This completes the proof in the case of an even $n$.\vv

 Let $n$ be odd. If  $s=1$ (i.e. the leftover is positive), then the case is reducible to the above case in virtue of Theorem \ref{fullor2}.2. So, let $s=-1$ (a negative leftover). Let us consider the chain of doubletons and integrate one of them with the leftover into a generator of length 3. Among three types of groups, $(3,0)\leftrightarrow (1,2)$ is of interest.
 By inspection, the integration/disintegration proceeds as follows:
 \[
 (++)(-)\,\longleftrightarrow\, (+++) \mbox{ (which is $(3,0)$), and}\quad
 (--)(-)\,\longleftrightarrow\, (--+) \mbox{ (which is $(1,2)$)}.
 \]
 Equivalently, in our symbolic notation:
 \[
 \PNs 1 0 - \,\longleftrightarrow\, (3,0)\leftrightarrow (1,2)\,\longleftrightarrow \,\PNs 0 1 -
 \]
 In other words, with $s=-1$ at hand we may switch the parity of any doubleton, still within the same group. This completes the proof in the odd case.
 \QED

\subsection{Dual decomposition}
It is natural to wonder about a dual decomposition, i.e. about the possibility of replacing $E$ by $E'=F_0\cup F$, where $F=F_1\cup\cdots \cup  F_k$ such that each $F_k$  is commutative for $k\ge 1$, $F_0$ is anticommutative, and $F_j\circ F_{j'}=-1$ for $j\neq j',\, 1\le j,\,j'\le k$.  $F_0$ would commute with other components.\vv

Conversely, when we face such pattern we may try to find an `explanation', i.e., a simple generator that `causes' it. We will show how simple counts help to fulfill this objective.
 \begin{example}
 Let $(\wek f,\wek g)$ have length $N$, where $\wek f\circ\wek g=-1$ and both sequences are nonempty and commutative. Then the count $c^-=3\cdot 2^{2N-3}$ and the pattern is `caused' by just two anticommuting elements.
 \end{example}
 \Proof Consider $\wek f^{\sp}\, \sg^{\sr}\circ \wek f^{\sq}\,\sg^{\wek s}=(-1)^{ps+qr}$. Let $\wek f$ have length $i$ and denote the length of $\sg$ by  $j=N-i$. Then
 \[
 c^-=\frac{1}{2}\sum_{\sp,\sr,\sq,\wek s} \Big(1-(-1)^{ps+qr}\Big)=\frac{2^{2N}-a^2}{2},
 \]
where the number $a$ is easily computable, yielding the claimed count:
\[
\begin{array}{rl}
a=\Sum_{\sq,\sr}(-1)^{qr}&=\Sum_{q\cong 0,\sr} 1+
\Sum_{q\cong 1,r \cong 0} 1
-\Sum_{q\cong 1, r\cong 1}1 \\
\vspandex &=2^{i-1}\cdot 2^j + 2^{i-1}\cdot 2^{j-1} -2^{i-1}\cdot 2^{j-1} =2^{N-1}.
\end{array}
\]
Now we look for the even $n$, the maximal size of an AC basic sequence $\wek e$ that together with $m=N-n$ commuting elements, also commuting with $\wek e$, forms a generator. In virtue of Theorem \ref{quickcounts} the AC count equals $2^{2m+n-1} \left(2^n-1\right)=2^{2N-n-1}\left(2^n-1\right)$. Since the quantity is a group invariant, it also equals $3 \cdot 2^{2N-3}$. For an even $n$, we easily check that $2^n-1$ is an odd multiple of 3. Hence, necessarily, $2N-n-1=2N-3$, so $n=2$. In other words, just a single pair of AC elements `causes' the given pattern, although not uniquely.
\QED
\vv
Several `recipes' follow. Let $e_1\circ e_2=-1$, $d_p\circ d_{q}=1$ for $p,q\in\set{1,...,N-2}$, and $\wek e\circ\wek d=1$.\vv

{\bf Recipe 1}.  In the first example we put
\[
f_1=e_1,\,f_{p}=e_1 d_p,\, p=2,\dots,i-1,\quad \mbox{and}\quad g_q=e_1e_2 d_q,\,q=N-j-1,...,N-2.
\]
Observe that $N-j-1=i-1$, i.e., the last member in the first group is $e_1d_{i-1}$ while the first member in the second group is $e_1e_2d_{i-1}$, so both generate $e_2$ and thus, together with $e_1$, all $d_p$'s. E.g.,  for $i=5$ and $j=3$:
\[
F=\left\{e_1,e_1d_1,e_1d_2,e_1d_3,e_1d_4\right\},\quad
G=\left\{e_1e_2d_4,e_1e_2d_5,e_1e_2d_6\right\}.
\]
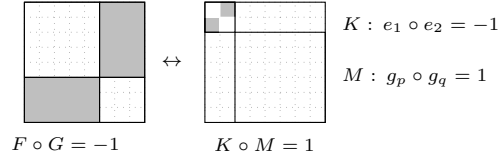
\begin{center}
\begin{tikzpicture}[scale=0.2]
\draw[help lines,dotted] (0,0) grid (8,8);
\draw [fill=lightgray] (5,3) rectangle (8,8);
\draw [fill=lightgray] (0,0) rectangle (5,3);
\draw[black,thin] (0,3) --(5,3) -- (5,8)  -- (0,8) -- (0,3);
\draw[black,thin] (5,0) --(8,0) -- (8,3)  -- (5,3) -- (5,0);
\draw[black,thin] (0,8) --(8,8) -- (8,0)  -- (0,0) -- (0,8);

\node[right] at (8.5,4) {\scriptsize $\leftrightarrow$};
\draw[help lines,dotted] (12,0) grid (20,8);
\draw [fill=lightgray] (12,6) rectangle (14,8);
\draw [fill=white,white] (12,7) rectangle (13,8);
\draw [fill=white,white] (13,6) rectangle (14,7);
\node[right] at (-1.5,-1.5) {\scriptsize $F\circ G=-1$};
\draw[black,thin] (12,6) --(14,6) -- (14,8)  -- (12,8) -- (12,6);
\draw[black,thin] (14,0) --(20,0) -- (20,6)  -- (14,6) -- (14,0);
\draw[black,thin] (12,8) --(20,8) -- (20,0)  -- (12,0) -- (12,8);
\node[right] at (20.5,6.5) {\scriptsize $K:\,e_1\circ e_2=-1$};
\node[right] at (20.5,3) {\scriptsize $M:\,g_p\circ g_q=1$};
\node[right] at (12,-1.5) {\scriptsize $K\circ M=1$};
\end{tikzpicture}
\end{center}
Here both sets are pure, $\sigma(F)=\sigma(e_1)$ while $\sigma(G)=-\sigma(e_1)\sigma(e_2)$.\vv

{\bf Recipe 2}. In another example the sets may have arbitrary signatures $\pm 1$. First, let us `pile' commuting elements as follows and then assign each product to one of two disjoint sets, $D_1$ of size $i$ or $D_2$ of size $j$, $i+j=N-2$:
\[
d_1,\,d_1d_2,\,\dots, d_1d_2\cdots d_{N-2}.
\]
Then the commuting sets $F=\set{e_1}\cup e_1D_1$ and $G=\set{e_2}\cup e_2D_2$ anticommute and are pure of signatures of $e_1$ or $e_2$. \vv

{\bf Recipe 3}. Alternatively, let us use $e_1$ and $e_2$ of opposite signatures to build an arbitrarily mixed $F$, and let $e_1e_2$ yield $G$, which must be pure.
\vv
{\bf General Recipe}. The latter recipe immediately generalizes to an arbitrary number of anticommuting commutative generators. The signatures can be controlled by signatures of commuting  elements $d_l$. Since the issue turns out to be rather elementary, we omit further details.
\vspace{10pt}

\noindent {\bf Acknowledgement.} Example \ref{2x2} emerged from a discussion with Dr.\ A.\ Jadczyk.

\end{document}